\documentclass[11pt]{article}
\usepackage{amsmath, amsthm,amsfonts,amssymb,amscd}
\usepackage{color}
\usepackage{geometry}
\usepackage{float}
\usepackage{soul}
\usepackage{marginnote}
\def\disp{\displaystyle}

\def\tto{\;{\lower 1pt \hbox{$\rightarrow$}}\kern -10pt
\hbox{\raise 2pt \hbox{$\rightarrow$}}\;}
\def\Hat{\widehat}
\def\hat{\widehat}

\def\tilde{\widetilde}
\def\Bar{\overline}
\def\ra{\rangle}
\def\la{\langle}
\def\ve{\varepsilon}
\def\B{\mathbb B}

\def\h{\hfill\Box}
\def\R{\mathbb R}

\def\ox{\bar{x}}

\def\oy{\bar{y}}

\def\ov{\bar{v}}

\def\gph{\mbox{\rm gph}\,}

\def\dom{\mbox{\rm dom}\,}

\def\h{\hfill\triangle}
\def\dn{\downarrow}
\def\O{\Omega}
\def\Lm{\Lambda}

\def\st{\stackrel}
\def\oR{\Bar{\R}}

\def\lm{\lambda}
\def\gg{\gamma}

\def\kk{\kappa}

\def\vt{\vartheta}
\newcounter{lk}

\begin{document}

\newtheorem{Theorem}{Theorem}[section]
\newtheorem{Proposition}[Theorem]{Proposition}
\newtheorem{Remark}[Theorem]{Remark}
\newtheorem{Lemma}[Theorem]{Lemma}
\newtheorem{Corollary}[Theorem]{Corollary}
\newtheorem{Definition}[Theorem]{Definition}
\newtheorem{Example}[Theorem]{Example}
\renewcommand{\theequation}{\thesection.\arabic{equation}}
\normalsize
\def\proof{
\normalfont
\medskip
{\noindent\itshape Proof.\hspace*{6pt}\ignorespaces}}
\def\endproof{$\h$ \vspace*{0.1in}}

\title{\small \bf CHARACTERIZATION OF TILT STABILITY VIA SUBGRADIENT GRAPHICAL DERIVATIVE WITH APPLICATIONS TO  NONLINEAR  PROGRAMMING}
\date{}
\author{Nguyen Huy  Chieu\footnote{Institute of Natural Sciences Education, Vinh
University,  Nghe  An, Vietnam; email: nghuychieu@gmail.com. }, \, \ Le Van  Hien\footnote{Department of Natural Science Teachers, Ha Tinh University,  Ha Tinh, Vietnam; email: lehiendhv@gmail.com.},
\, \  Tran T.A. Nghia\footnote{Department of Mathematics and Statistics, Oakland University, Rochester, MI 48309, USA; email: nttran@oakland.edu.}}
\maketitle
{\small \begin{abstract}
This paper  is  devoted to the study of tilt stability   in finite dimensional optimization via the approach of using the subgradient graphical derivative.  We    establish a new characterization of  tilt-stable local minimizers  for a broad class of  unconstrained optimization problems  in terms of   a uniform  positive definiteness of  the subgradient  graphical derivative of the objective function around the  point in question.  By applying this result to  nonlinear programming under the metric subregularity constraint qualification,  we    derive a  second-order characterization and several new sufficient conditions for tilt stability. In particular, we show that each stationary point of a nonlinear programming problem satisfying the metric subregularity constraint qualification  is a tilt-stable local minimizer if  the classical strong second-order sufficient condition holds.
\end{abstract}}

\medskip

\normalsize
\section{Introduction}
\setcounter{equation}{0}

{\it Tilt stability}  is a property of local minimizers guaranteeing the minimizing point shifts in a Lipschitzian manner under  linear perturbations on the objective function of a optimization problem, which is a desired behavior in optimization from both theoretical and numerical viewpoints.
This notion was  introduced by Poliquin and Rockafellar~\cite{PR2} for problems of unconstrained optimization with extended-real-valued objective function. As usual, incorporating constraints into the objective function via the indicator function of the feasible set, one can speak of tilt stability for constrained optimization problems. Tilt stability is basically equivalent to  uniform second-order growth condition as well as strong metric regularity of the subdifferential \cite{BS, DL, MN}. These properties have been intensively studied  in the recent years; see  \cite{DL,  EW, MG15, LPR, LZ, MN, MO1, MR}.

The first characterization of tilt-stability using  second-order generalization differentiation was due to  Poliquin and Rockafellar~\cite{PR2}. In that paper they  proved that for an unconstrained optimization problem, under mild assumptions of prox-regularity,  a stationary point is a tilt-stable local minimizer if and only if  the second-order limiting subdifferential/generalized Hessian  in the sense of Mordukhovich \cite{M92}  is positive-definite at the point in question. Furthermore, using this result together with a formula of Dontchev and Rockafellar \cite{DR1} for the second-order limiting subdifferential of  the indicator function of polyhedral convex set, they obtained  a second-order characterization of tilt-stability for nonlinear programming problems with linear constraints~\cite[Theorem 4.5]{PR2}. The main difficulty in applying the tilt-stability characterization of Poliquin and Rockafellar~\cite{PR2}   to other nonlinear constrained optimization problems is the computation/estimation of the  second-order subdifferential in terms of explicit problem data.

By establishing  new second-order subdifferential calculi, Mordukhovich and Rockafellar~\cite{MR} derived second-order characterizations of tilt-stable minimizers for some  classes of constrained optimization problems. Among other important things, they showed that for  $C^2$-smooth  nonlinear programming problems, under the {\em linear independence constraint qualification} (LICQ), a stationary point is a  tilt-stable local minimizer if and only if the {\em strong second-order sufficient condition} (SSOSC) holds.  Consequently, in this setting, tilt-stability is equivalent to Robinson's strong regularity \cite{Ro} of the associated Karush-Kuhn-Tucker system whenever LICQ occurs at the point in question. In contrast to Robinson's strong regularity, tilt stability does not postulate   LICQ as a necessary condition. This observation motivated the study of tilt-stability for nonlinear programming under constraint qualifications weaker than LICQ, aiming to   cover a broader class of examined problems.

    Under  the validity of both the {\em Mangasarian-Fromovitz constraint qualification} (MFCQ) and the {\em constant rank constraint qualification} (CRCQ),   Mordukhovich and  Outrata \cite{MO1} proved  that SSOSC is a sufficient condition  for a stationary point to be a tilt-stable local minimizer in nonlinear programming. In \cite{MN} Mordukhovich and Nghia showed that SSOSC is indeed not a necessary condition for tilt stability and then introduced the {\em uniform second-order sufficient condition} (USOSC) to characterize tilt stability when both MFCQ and CRCQ occur. Recently,  Gfrerer and Mordukhovich~\cite{MG15} obtained some  pointbased second-order sufficient conditions for tilt-stable local minimizers under some  weak conditions including the so-called {\em metric subregularity constraint qualification} (MSCQ) and the
{\em bounded extreme point property} (BEPP).   Furthermore, some  pointbased second-order characterizations of tilt-stability were established in~\cite{MG15}  under certain additional assumptions.  We also note that  the  approach of using a uniform positive definiteness of  the {\em combined second-order subdifferential} by Mordukhovich and Nghia~\cite{MN} is an essential tool for the analysis done in~\cite{MG15}. For more information on the recent literature on tilt stability in nonlinear programming, we refer the reader to  \cite{MG15, MN,  MO1, MS16} and the references therein.

The approach of this paper is different from those in the aforementioned references. We mainly use the {\it subgradient graphical derivative} of an  extended-real-valued function, which is  the {\em graphical derivative} of its  limiting subdifferential \cite{rw},  to characterize tilt stability and extend to the case of nonlinear programming.  In fact, the subgradient graphical derivative and the second-order subdifferential are generally independent concepts; however, under additional conditions, the value of the subgradient graphical derivative can be identified with a subset of the value of the second-order subdifferential; see \cite{rw, rz}. We note that one of the biggest advantages of this approach is the workable computation of  the graphical derivative in  various important cases under very mild assumptions in initial data; see \cite{CH16, GO, GY16, MOR1}. Furthermore,  several results on tilt stability, e.g.  in \cite{MG15},  were established based on the calculation of   the subgradient graphical derivative as a mediate step.  These observations  indeed lead us to the following natural questions:

{\it Is it possible to use the subgradient graphical derivative to characterize  tilt stability of local minimizers for unconstrained optimization problem? If yes, is such a characterization  useful in helping us to improve the knowledge of tilt stability  for nonlinear programming problems?}

The aim of this paper is to give the positive answers for  the two risen questions. Precisely, after recalling some preliminary materials in Section 2,  we   establish   a  new second-order characterization of  tilt-stable local minimizers  for  unconstrained optimization problems, in which the objective function is prox-regular and subdifferentially continuous \cite{PR2} in Section~3. The characterization is expressed in terms of   a uniform  positive definiteness
 of  the subgradient graphical derivative of the objective function around the considered point in which  the prox-regularity of the objective function is essential not only for the necessary implication but also for the sufficient one.  In Section 4, by  applying the established characterization to  nonlinear programs under the metric subregularity constraint qualification,  we  derive a  new second-order characterization of tilt stability via a  uniform second-order sufficient condition, which reduces to \cite[Theorem 4.3]{MN} under the validity of both MFCQ and CRCQ
  and then obtain pointbased second-order sufficient conditions  for  a stationary point of the problem  to be a tilt-stable local minimizer. As a consequence, we show that  each stationary point of a nonlinear programming problem satisfying  MSCQ is a tilt-stable local minimizer if  SSOSC is satisfied. This result improves the corresponding result of Mordukhovich and  Outrata \cite{MO1} replacing the combination of MFCQ and CRCQ by the much weaker MSCQ. The final Section 5 involves some perspectives of the obtained results and future works.

\section{Preliminaries}
\setcounter{equation}{0}

In this section we recall some basic notions and  facts  from variational analysis that will be used repeatedly  in the sequel; see  \cite{DR14,M1, rw} for more details. Let $\O$ be a nonempty subset of the Euclidean space $\R^n$ and $\ox$ be a point in $\O$.  Define the polar cone of $\O$  by  $\Omega^\circ:=\big\{v\in \R^n|\, \langle v,x\rangle\leq 0\ \, \mbox{for all}\ x\in \Omega\big\}.$
The (Bouligand-Severi) {\it tangent/contingent cone} to the set $\Omega$ at $\bar x\in \Omega$  is  known as
$$T_\Omega(\bar x):=\big\{v\in\R^n | \, \mbox{there exist}\  t_k \downarrow 0, \  v_k\rightarrow v\ \mbox{ with }\  \bar x+t_kv_k\in\Omega \ \mbox{for all} \   k\in \mathbb{N}\big\}.$$
The  polar cone  of the tangent cone  is  the  (Fr\'echet) {\it regular normal cone} to $\Omega$ at $\bar x \in \Omega$ defined by
 \begin{equation}\label{eqdualTN}\widehat{N}_\Omega(\bar x):=T_\Omega(\bar x)^\circ.\end{equation}
It is well-known that the regular normal cone could be presented by the following construction
$$\widehat{N}_\Omega(\bar x):=\left\{v\in\R^n\, \big|\, \limsup\limits_{x\st{\Omega}\rightarrow\bar x}\frac{\langle v, x-\bar x\rangle}{\|x-\bar x\|}\leq 0\right\},$$
where  $x\st{\Omega}\rightarrow\bar x$ means that $x\rightarrow \bar x$ with $x\in \Omega.$
Another normal cone construction used in our work is the  (Mordukhovich) {\it limiting/basic  normal cone} to $\Omega$ at $\bar x\in \Omega$ defined by
$$N_\Omega(\bar x)=\big\{v\in \R^n\, | \,  \mbox{there exist}\,   x_k\st{\Omega}\rightarrow \bar x,\, v_k\in \widehat{N}_\Omega(x_k) \  \mbox{with} \  v_k\rightarrow v\big\},$$
which was introduced   by Mordukhovich \cite{M76} in an equivalent form.
If $\bar x\not\in \Omega,$ one puts $T_\Omega(\bar x)=\emptyset$ and $N_\Omega(\bar x)=\widehat{N}_\Omega(\bar x)=\emptyset$ by convention.
When the set $\Omega$ is convex, the above tangent cone and normal cones  reduce to the tangent cone and normal cone in the  sense of classical convex analysis.

Consider the set-valued mapping $F: \R^n\rightrightarrows\R^m$ with the domain $\dom F:=\big\{x\in \R^n|\; F(x)\neq \emptyset\big\}$ and graph ${\rm gph}\, F:=\big\{(x,y)\in \R^n\times \R^m  |\,   y\in F(x)\big\}$. Suppose that the domain of $F$ is nonempty and $(\ox,\oy)$ is an element of ${\rm gph}\, F$.

The {\it graphical derivative} of $F$ at $\bar x$ for $\bar y\in F(\bar x)$ is the set-valued mapping $DF(\bar x|\bar y): \R^n\rightrightarrows\R^m$ defined by
\begin{equation*}\label{GraDer}
DF(\bar x|\bar y)(w):=\big\{ z\in\R^m \,  |\,  (w,z)\in T_{{\rm gph}F}(\bar x,\bar y)\big\} \ \,\, \mbox{for }\   w\in\R^n,
\end{equation*}
  that is,  ${\rm gph}\, DF(\bar x|\bar y)= T_{{\rm gph}F}(\bar x,\bar y).$  This concept was introduced in the early  1980s by Jean-Pierre Aubin, who called it the contingent derivative. Here we follow the references \cite{DR14, rw} in using the terminology ``the graphical derivative''.    In the case  $F(\bar x)=\{\bar y\},$  one writes  $DF(\bar x)$   for  $DF(\bar x|\bar y)$.
  If $(\bar x,\bar y)\not\in \gph F,$ one puts  $DF(\bar x|\bar y)(w)=\emptyset$ for all $w$ by convention.
   We note further that if $\Phi: \R^n\rightarrow \R^m$ is a single-valued mapping differentiable at $\bar x,$  then $ D\Phi(\bar x)=\nabla \Phi(\bar x).$


Recall \cite{M1, rw} that  the set-valued  mapping $F:\R^n\tto \R^m$ is said to be Lipschitz-like (pseudo-Lipschitz or has the Aubin property) at $\ox\in {\rm dom}\, F$ for $\oy\in F(\ox)$ with modulus $\ell\ge 0$ if there exist neighborhoods (nbhs) $U$ of $\ox$ and $V$ of $\oy$ such that
\begin{equation}\label{LiLi}
F(x)\cap V\subset F(u)+\ell\|x-u\|\B\quad \mbox{for all } \quad x,u\in U,
\end{equation}
where $\B$ is the unit ball in $\R^m$.  The exact modulus for Lipschitz-like property of $F$ at $\ox$ for $\oy$ is defined by
\begin{equation*}
{\rm lip}\, F(\ox|\oy):=\inf\big\{\ell\in \R_+|\; \exists \mbox{ nbhs } U \mbox{ of } \ox \mbox{ and } V \mbox{ of } \oy \mbox { such that \eqref{LiLi} holds}\big\}. \label{elip}
\end{equation*}
It is known from \cite[Theorem 4B.2]{DR14} that $F$ is Lipschitz-like at $\ox$ for $\oy$ if and only if
\begin{equation}\label{Don}
\limsup_{(x,y)\st{{\rm gph F}}\to(\ox,\oy)}|DF(x|y)|^{-}<\infty,
\end{equation}
where $\disp|DF(x|y)|^{-}=\sup_{\|w\|\le 1}\inf_{z\in DF(x|y)(w)}\|z\|$. Moreover, the quantity on the left-hand side of~\eqref{Don} is the exact modulus of $F$ at $\ox$ for $\oy$.

An important property of  set-valued mapping known as {\em metric regularity} also plays essential roles in our study. The set-valued  mapping $F$ is said to be metrically regular at $\ox\in {\rm dom}\, F$ for $\oy\in F(\ox)$ with modulus $\kk> 0$ if there exist nbhs $U$ of $\ox$ and $V$ of $\oy$ such that
\begin{equation}\label{met-reg}
d(x;F^{-1}(y))\le\kk d(y;F(x))\quad \mbox{for all}\quad  (x,y)\in U\times V,
\end{equation}
where $d(x;\O)$ represents the distance from a point $x\in \R^n$ to a set $\O\subset \R^n$.  The infimum of all such $\kk$ is the modulus of metric regularity, denoted by ${\rm reg}\, F(\ox|\oy)$.
 It is well-known that $F$ is metrically regular at $\ox$ for $\oy$ with modulus $\kk>0$ if and only if $F^{-1}$ is Lipschitz-like at $\oy$ for $\ox$ with the same modulus; see, e.g., \cite[Theorem~1.49]{M1}.

Following \cite[Section 3.8]{DR14}, we say $F$ is   {\em metrically  subregular} at $\ox\in {\rm dom}\, F$ for $\oy\in F(\ox)$ with modulus $\kk> 0$ when the inequality \eqref{met-reg} holds with $y=\oy$, i.e.,
\begin{equation*}\label{subreg}
d(x;F^{-1}(\oy))\le\kk d(\oy;F(x))\quad \mbox{for all}\quad  x\in U.
\end{equation*}
The infimum of all such $\kk$ is the modulus of metric subregularity, denoted by ${\rm subreg}\, F(\ox|\oy)$.

The set-valued mapping $F$ is said to be {\em strongly metrically regular} at $\ox$  for $\oy\in F(\ox)$ with modulus $\kk>0$ if its inverse $F^{-1}$ admits a {\em single valued and Lipschitz continuous localization} around $\ox$ for $\oy$ with modulus $\kk> 0$ in the following sense:  there are neighborhoods $U$ of $\ox$ and $V$ of $\oy$ and a Lipschitz continuous function  $\vartheta:V\to U$ with full domain $U$ and constant $\kk$ satisfying that
\[
\gph \vartheta=\gph F^{-1}\cap (U\times V).
\]
Strong metric regularity introduced by Robinson \cite{Ro} has been known a strong notion useful in optimization and algorithm; see \cite{DR14} for further discussions and applications to nonlinear programming.

Assume that  $f: \R^n\rightarrow \overline{\R}:=\R\cup\{\infty\}$ is an extended-real-valued proper function with $\bar x\in\dom f:=\big\{x\in \R^n|\; f(x)<\infty\big\}$. The {\it limiting  subdifferential} (known also as the Mordukhovich/basic subdifferential)  of $f$ at $\bar x$ is defined by
$$\partial f(\bar x):=\big\{ v\in\R^n\, |\, (v,-1)\in N_{{\rm epi} f}(\bar x,\bar y)\big\},$$
where ${\rm epi}\, f:=\big\{(x,r)\in \R^n\times \R|\;r\ge \varphi(x)\big\}$ is the epigraph of $f$.
 For each $(x,v)\in \R^n\times \R^n,$  following \cite{MOR},  the mapping $D\partial f(x,v): \R^n\rightrightarrows \R^n$ defined by
 \[D\partial f(x,v)(w):=\big\{z\, |\, (w,z)\in T_{\gph\partial f}(x,v)\big\}
 \]  is said to be the {\it subgradient graphical derivative} of $f$ at $x$ for $v.$
 Finally let us recall two notations of prox-regularity and subdifferential continuity. Function $f$ is  said to be {\it prox-regular} at $\ox\in \dom f$ for $\ov\in \partial f(\ox)$ if there exist $r,\ve>0$ such that for all $x,u\in \B_\ve(\ox)$ with $|f(u)-f(\ox)|<\ve$ we have
\begin{equation}\label{prox}
f(x)\ge f(u)+\la v, x-u\ra-\frac{r}{2}\|x-u\|^2\quad \mbox{for all}\quad v\in \partial f(x)\cap \B_\ve(\ov).
\end{equation}
Moreover, we say $f$ is {\em subdifferentially continuous} at $\ox$ for $\ov$ if the mapping $(x,v)\mapsto f(x)$ is continuous relative to the graph of $\partial f$ at $(\ox,\ov)$. When the function $f$ is both prox-regular and subdifferentially continuous at $\ox$ for $\ov$, by choosing smaller $\ve>0$, \eqref{prox} is still valid without the restriction ``$|f(u)-f(\ox)|<\ve$''. It is worth noting that in this case the graph of $\partial f$ is closed around $(\ox,\ov)$. For more detailed information on the prox-regularity and its applications, we refer the reader to the references \cite{CT10, PR1, rw}.

 \section{Second-Order Characterizations of Tilt Stability}
\setcounter{equation}{0}

This section focuses on the tilt stability for unconstrained optimization problems.  The  concept of tilt-stability due to  Poliquin and Rockafellar \cite{PR2} is defined  as follows.

\begin{Definition}{\bf (Tilt stability \cite{PR2}).}\label{tilt}  Given $f\colon\R^n\to\oR$, a point $\ox\in\dom f$ is a {\sc tilt-stable local minimizer} of $f$ with modulus $\kk>0$ if there is a number $\gg>0$ such that the mapping
\begin{eqnarray*}\label{2.10}
M_\gg:v\mapsto{\rm argmin}\big\{f(x)-\la v,x\ra\big|\;x\in\overline\B_\gg(\ox)\big\}
\end{eqnarray*}
is single-valued and Lipschitz continuous with constant $\kk>0$ on some neighborhood of $0\in \R^n$ with $M_\gg(0)=\ox$. In this case we define the exact modulus for tilt stability of function $f$ at $\ox$ by
\[
{\rm tilt}\, (f,\ox):=\inf\big\{\kk|\; \ox \mbox{ is a tilt-stable minimizer of $f$ with modulus $\kk>0$}\big\}.
\]
\end{Definition}

The following result taken from \cite[Theorem 3.1 and Theorem~3.2]{MN} provides some useful characterizations for tilt stability via the strong metric regularity of the subdifferential and the uniform second-order growth condition; see also \cite[Theorem~3.3]{DL} for the earlier result without paying much attention to the modulus of tilt-stability.

\begin{Theorem}{\bf (Characterizations of tilt stability).}\label{thm0} Let $f:\R^n\to\oR$ be a lower semi-continuous (l.s.c.) proper function such that  $\ox\in \dom f$ and  $0\in\partial f(\ox)$. Assume that $f$ is both prox-regular and subdifferentially continuous at $\ox$ for $\ov=0$. Then the following assertions are equivalent:

{\bf (i)} The point $\ox$ is a tilt-stable local minimizer of the function $f$ with modulus $\kk>0$.

{\bf (ii)} The point $\ox$ is a  local minimizer of $f$ and $\partial f$ is {\sc strongly metrically regular} at $\ox$ for $\ov$ with modulus $\kk>0$ in the sense that $(\partial f)^{-1}$ admits a single-valued and Lipschitz continuous localization around $\ov$ for $\ox$ with modulus $\kk>0$.

{\bf (iii)} There are neighborhoods $U$ of $\ox$ and $V$ of $\ov$ such that the mapping $(\partial f)^{-1}$ admits a single-valued localization $\vt:V\to U$ around $\ov$ for $\ox$ and that for any pair $(v,u)\in \gph \vt=\gph(\partial f)^{-1}\cap (V\times U)$ we have the {\sc uniform second-order growth condition}
\begin{eqnarray}\label{3.1a}
f(x)\ge f(u)+\la v,x-u\ra+\frac{1}{2\kk}\|x-u\|^2\;\mbox{ whenever }\;x\in U.
\end{eqnarray}
\end{Theorem}

Tilt stability has been  also characterized  via second-order subdifferentials,  in particular,  the limiting second-order subdifferential, that is, the limiting coderivative to the limiting subdifferential; see, e.g., \cite{EW,  MN, MG15, MOR, MR, MS16,   PR2}. Over the years,  this dual approach has produced many nice results on tilt stability,  leading  to various  applications to nonlinear programming, semidefinite programming, conic programming and so on. To the best of our knowledge, in the current stage, the dual approach has met some severe difficulties in handling  tilt stability for  non-polyhedral conic programs under weak conditions, due to the limitation of  computing such dual second-order structures under mild assumptions.

We next  examine  a new  approach  to tilt stability, which is  based on  {\it the subgradient graphical derivative}.  It turns out that, as shown  in the next section, this approach can help us to improve the knowledge of tilt stability for nonlinear programming problems. Precisely, we have the following  theorem  which provides a new second-order characterization of tilt stability that will be the main tool in investigating tilt stability for nonlinear programming problems in Section~\ref{TSNP}.

\begin{Theorem}{\bf (Subgradient graphical derivative characterization of tilt-stability).}\label{thm1}
Let $f:\R^n\to\oR$ be a l.s.c.\ proper function with $\ox\in \dom f$ and  $0\in\partial f(\ox)$. Assume that $f$ is both prox-regular and subdifferentially continuous at $\ox$ for $\ov=0$. Then the following assertions are equivalent:

{\bf (i)} The point $\ox$ is a tilt-stable local minimizer of $f$ with modulus $\kk>0$.

{\bf (ii)} There is a constant $\eta>0$ such that for all $w\in\R^n$ we have
\begin{eqnarray}\label{3.2}
\la z,w\ra\ge\frac{1}{\kk}\|w\|^2\;\mbox{ whenever }\;z\in D\partial f(u,v)(w)\;\mbox{ with }\;(u,v)\in\gph\partial f\cap\B_\eta(\ox,0).
\end{eqnarray}
Furthermore, the exact tilt-stable modulus of $f$ is calculated by the formula
\begin{eqnarray}\label{3.2b}
{\rm tilt}\, (f,\ox)=\inf_{\eta>0}\sup\Big\{\frac{\|w\|^2}{\la z,w\ra}\Big|\; z\in D \partial f (u|v)(w),\;(u,v)\in\gph\partial f\cap\B_\eta(\ox,0)\Big\}
\end{eqnarray}
with the convention that $0/0=0$.
\end{Theorem}
\noindent{\bf Proof.} To verify {\bf (i)}$\Longrightarrow${\bf (ii)}, suppose that $\ox$ is a tilt-stable local minimizer of $f$ with modulus $\kk>0$. Then we get from Theorem~\ref{thm0} that there exists  a single-valued and Lipschitz continuous  localization $\vt$ of $(\partial f)^{-1}$ relative to a neighborhood $V\times U$ of $(\ov,\ox)$ such that \eqref{3.1a} is satisfied. Fix  $\gg>0$ with $\B_\gg(\ox)\subset U$, due to the Lipschitz continuity of $\vt$ and $\vt(\ov)=\ox$ we find some $\nu>0$ such that $\vt({\rm int\,}\B_\nu(\ov))\subset {\rm int\,}\B_\gg(\ox)$.  Since $\ox$ is a tilt-stable local minimizer of the function $f$ with modulus $\kk>0$, the positive real constants  $\gg$ and   $\nu$ can be chosen such that $M_\gg$ is single-valued and  Lipschitz continuous with modulus $\kk$ over $\B_\nu(\ov).$
It follows from~\eqref{3.1a} that
   $$f(x)-\la v, x\ra \geq  f\big(\vt(v)\big)-\la v, \vt(v)\ra\ \, \mbox{for all}\ x\in \B_\gg(\ox)\ \mbox{and}\ v\in\B_\nu(\ov).$$
Then we have $\vt(v)\in  M_\gg(v)$ and thus $M_\gg(v)=\vt(v)$ for all $v\in {\rm int\,}\B_\nu(\ov)$.  The latter is due to the single-valuedness of $M_\gg$ over $\B_\nu(\ov).$
Set $\tilde U:={\rm int\,}\B_\gg(\ox)$ and $\tilde V:={\rm int\,}\B_\nu(\ov)$. Take any $(x,y),\,(u,v)\in\gph\partial f\cap(\tilde U\times \tilde V)$, we get from \eqref{3.1a} that
\begin{eqnarray*}\begin{array}{ll}
f(x)&\ge f(u)+\la v,x-u\ra+\frac{1}{2\kk}\|x-u\|^2\\
f(u)&\ge f(x)+\la y,u-x\ra+\frac{1}{2\kk}\|u-x\|^2.
\end{array}
\end{eqnarray*}
Adding these two inequalities gives  us that
\begin{eqnarray}\label{3.6}
\la y-v,x-u\ra\ge \frac{1}{\kk}\|x-u\|^2\;\mbox{ whenever  }\;(x,y),\,(u,v)\in\gph\partial f\cap(\tilde U\times\tilde  V).
\end{eqnarray}
To justify \eqref{3.2}, pick any $z\in D\partial f(u|v)(w)$ with $(u,v)\in\gph\partial f\cap (\tilde U\times \tilde V)$ and $w\in\R^n$ and find sequences $t_k\dn 0$ and $(w_k,z_k)\to (w,z)$ such that $(u,v)+t_k(w_k,z_k)\in \gph  \partial f\cap (\tilde U\times \tilde V).$  Hence, we derive from \eqref{3.6} that
\[
t_k^2\la z_k,w_k\ra=\la v+t_kz_k-v,u+t_kw_k-u\ra\ge \frac{1}{\kk}\|u+t_kw_k-u\|^2=\frac{1}{\kk}t_k^2\|w_k\|^2
\]
for all  $k$. It follows that $\la z_k,w_k\ra\ge \frac{1}{\kk} \|w_k\|^2$. Taking $k\to \infty$ in the latter inequality gives us that $\la z,w\ra\ge \frac{1}{\kk}\|w\|^2$, which clearly ensures \eqref{3.2} with $\eta=\min\{\gg,\nu\}$. \vspace*{0.05in}

Let us now justify the converse implication {\bf (ii)}$\Longrightarrow${\bf (i)} by supposing that \eqref{3.2} holds with some $\eta,\kk>0$. Since $f$ is prox-regular and subdifferentially continuous at $\ox$ for $\ov$, there are $r,\ve>0$ with $\ve<\eta$ satisfying
\begin{eqnarray}\label{3.7a}
f(x)\ge f(u)+\la v,x-u\ra-\frac{r}{2}\|x-u\|^2,
\end{eqnarray}
for all $x\in\B_\ve(\ox)$ and  $(u,v)\in\gph\partial f\cap\B_\ve(\ox,\ov).$ Pick any $s>r$ and define
\begin{equation}\label{Gi}
g(x):=f(x)+\frac{s}{2}\|x-\ox\|^2\quad\mbox{for all}\quad x\in\R^n.
\end{equation}
We have
\begin{equation}\label{sGi}\partial g(x)=\partial f(x)+s(x-\ox).
\end{equation}
Define further $W:=J(\B_\ve(\ox,\ov))$ with $J(u,v):=(u,v+s(u-\ox))$ for $(u,v)\in\R^n\times\R^n$ and observe that $W$ contains an open ball $\B_{2\delta}(\ox,\ov)$ for some $\delta\in (0,\ve) $ sufficiently small.
Take any $x\in\B_\delta(\ox)$ and $(u,v)\in\gph\partial g\cap\B_\delta(\ox,\ov).$ Then $\big(u, v-s(u-\ox)\big)=J^{-1}(u,v)\in \gph\partial f\cap\B_\ve(\ox,\ov).$ So by \eqref{3.7a} we have
$$f(x)\ge f(u)+\la v-s(u-\ox) ,x-u\ra-\frac{r}{2}\|x-u\|^2. $$
This together with \eqref{Gi}  implies  that
\begin{eqnarray}\label{3.8}
g(x)\ge g(u)+\la v,x-u\ra+\frac{s-r}{2}\|x-u\|^2
\end{eqnarray}
whenever  $x\in\B_\delta(\ox)$ and $(u,v)\in\gph\partial g\cap\B_\delta(\ox,\ov).$
Since $f$ is both prox-regular and subdifferentially continuous at $\ox$ for $\ov=0$,  it is easy to check from definition that $g$ is also prox-regular and subdifferentially continuous at $\ox$ for $\ov=0$.   Take any $(u,v)\in \gph\partial g\cap\B_\delta(\ox,\ov)$ and $z\in D\partial g(u|v)(w).$ By \cite[Proposition 4A.2]{DR14}, we get from \eqref{sGi} that
\[
z-sw\in D\partial f\big(u|v-s(u-\ox)\big)(w).
\]
Note that $(u,v-s(u-\ox))=J^{-1}(u,v)\subset J^{-1}(W)=\B_\ve(\ox,\ov)\subset\B_\eta(\ox,\ov)$. The above inclusion together with \eqref{3.2} implies that $\la z-sw,w\ra\ge\kk^{-1}\|w\|^2$. Thus
\begin{eqnarray*}\label{3.9}
\|z\|\cdot\|w\|\ge\la z,w\ra\ge(s+\kk^{-1})\|w\|^2,
\end{eqnarray*}
which ensures that $\|z\|\ge(s+\kk^{-1})\|w\|$. By the regular graphical derivative criterion for Lipschitz-like property from \cite[Theorem~4B.2]{DR14}  (see also \eqref{Don}) we conclude that $(\partial g)^{-1}$ is Lipschitz-like at $\ov$ for $\ox$ with the modulus $\ell:=(s+\kk^{-1})^{-1}$. Hence, there exists some neighborhood $U\times V\subset\B_\delta(\ox,\ov)$ of $(\ox,\ov)$ such that
\begin{equation}\label{LL}
(\partial g)^{-1}(v_1)\cap U\subset (\partial g)^{-1}(v_2) +\ell\|v_1-v_2\|\B\quad \mbox{for all}\quad v_1,v_2\in V.
\end{equation}
Since $\ox\in (\partial g)^{-1}(\ov)\cap U$, we get from the latter that
\begin{equation*}
\ox\in (\partial g)^{-1}(v)+\ell\|v-\ov\|\B\quad \mbox{for all}\quad v\in V.
\end{equation*}
So,  choosing the neighborhood $V$ smaller if necessary, assume that   $(\partial g)^{-1}(v)\cap \B_\beta(\bar x)\not=\emptyset$ for all $v\in V$ and  $\ell.{\rm diam}V<\beta,$
 where $\beta$ is some positive real number satisfying $\B_{2\beta}(\bar x)\subset U.$
Let $T$ be a localization of $(\partial g)^{-1}$ relative to $V\times U$, i.e., $\gph T=\gph (\partial g)^{-1}\cap(V\times U)$.  Since $\emptyset\neq(\partial g)^{-1}(v)\cap \B_\beta(\bar x)\subset (\partial g)^{-1}(v)\cap U$ for  $v\in V$, we have $\dom T=V$.   Moreover,  for all $(v_1,u_1), (v_2,u_2)\in\gph T=\gph(\partial g)^{-1}\cap (U\times V) \subset \gph(\partial g)^{-1} \cap\B_\delta(\ov,\ox)$,  it  follows from \eqref{3.8} that
\begin{eqnarray*}\begin{array}{ll}
\|v_2-v_1\|\cdot \disp\|u_2-u_1\|&\geq \la v_2-v_1, u_2-u_1\ra=\la v_2,u_2-u_1\ra+\la v_1, u_1-u_2\ra\\
&\disp\geq g(u_2)-g(u_1)+ \frac{s-r}{2}\|u_1-u_2\|^2+g(u_1)-g(u_2)+ \frac{s-r}{2}\|u_2-u_1\|^2\\
&\ge(s-r)\|u_2-u_1\|^2.
\end{array}
\end{eqnarray*}
This easily implies  that $T$ is single-valued, i.e., $v_k=T(u_k)$ for $k=1,2$.
 Since $(\partial g)^{-1}(v)\cap \B_\beta(\bar x)\not=\emptyset,$ we have $\emptyset\neq T(v)\in \B_\beta(\bar x)\subset U$ for all $v\in V.$ Furthermore, from  \eqref{LL} we obtain that $$T(v_1)\in (\partial g)^{-1}(v_2) +\ell\|v_1-v_2\|\B,$$or equivalently,  there exists $x_2\in (\partial g)^{-1}(v_2) $ such that
 \begin{equation}\label{seLi}
 \|T(v_1)-x_2\|\leq \ell \|v_1-v_2\|.
 \end{equation}
Hence,$$\begin{array}{rl}\|x_2-\ox\|&\leq \|x_2-T(v_1)\|+\|T(v_1)-\ox\|\\&\leq \ell \|v_1-v_2\|+\beta\\& \leq  \ell.{\rm diam}V +\beta<2\beta.\end{array}$$
It implies  that$$x_2\in  (\partial g)^{-1}(v_2)\cap \B_{2\beta}(\bar x)\subset   (\partial g)^{-1}(v_2)\cap U=T(v_2),$$
which together with \eqref{seLi} ensures that
\begin{eqnarray*}
\|T(v_1)-T(v_2)\|\leq \ell\|v_1-v_2\|\ \, \mbox{for all}\ v_1,v_2\in V.
\end{eqnarray*}
Taking into account that $g$ satisfies \eqref{3.8}, we get from Theorem~\ref{thm0} the existence of an open neighborhood $\hat U\times \hat  V$ of $(\ox,\ov)$ such that
\begin{eqnarray*}\label{3.12a}
g(x)\ge g(u)+\la v,x-u\ra+\frac{1}{2\ell}\|x-u\|^2\;\mbox{ if }\;x\in \hat  U,\,(u,v)\in\gph\partial g\cap(\hat  U\times \hat  V).
\end{eqnarray*}
Since $J^{-1}(\hat  U\times \hat  V)$ is an open set of $(\ox,\ov)$, we find a neighborhood $\tilde  U\times \tilde V$ of $(\ox,\ov)$ with $J(\tilde  U\times \tilde V)\subset \hat  U\times \hat  V$ and $\tilde U\subset \hat  U$. For any  $(u,v)\in\gph\partial f\cap  (\tilde U\times \tilde V)$ and $x\in \tilde U\subset \hat  U$, we have
\begin{eqnarray*}\begin{array}{ll}
f(x)&\disp=g(x)-\frac{s}{2}\|x-\ox\|^2\ge g(u)+\la v+s(u-\ox),x-u\ra+\frac{s+\kk^{-1}}{2}\|x-u\|^2-\frac{s}{2}\|x-\ox\|^2\nonumber\\
&=\disp g(u)+\la v,x-u\ra+\frac{1}{2\kk}\|x-u\|^2+\la s(u-\ox),x-u\ra+\frac{s}{2}\|x-u\|^2-\frac{s}{2}\|x-\ox\|^2\\
&\disp= g(u)+\la v,x-u\ra+\frac{1}{2\kk}\|x-u\|^2-\frac{s}{2}\|u-\ox\|^2\\
&=\disp f(u)+\la v,x-u\ra+\frac{1}{2\kk}\|x-u\|^2.
\end{array}
\end{eqnarray*}
Applying Theorem~\ref{thm0} again verifies {\bf (i)}. It is easy to see  that the exact bound formula~\eqref{3.2b} follows directly from \eqref{3.2}.\endproof

The next two examples show  the prox-regularity assumption  is essential  not only for   ${\bf (i)}\Rightarrow {\bf(ii)}$ but also for  ${\bf (ii)}\Rightarrow {\bf(i)}$ in  our Theorem~\ref{thm1}.

 \begin{Example}{\rm {\bf (Implication (i) $\Rightarrow$ (ii) fails in the absence of the prox-regularity).}\label{exam3.4}
Let $f: \R\rightarrow \R$ be the function defined by
\begin{equation*}  f(x) := \begin{cases}\   \min\bigg\{\big(1+\dfrac{1}{n}\big)|x|-\dfrac{1}{n(n+1)}, \dfrac{1}{n}\bigg\} \quad \  \mbox{if} \,      \  \dfrac{1}{n+1}\leq|x|\leq\dfrac{1}{n},\ n\in\mathbb{N^\ast},\\
 \quad  \ 0       \quad\quad\quad\quad\quad\quad\quad \quad  \ \ \ \ \ \ \ \ \ \ \ \ \ \ \ \ \ \ \mbox{if}    \,   \  x=0,\\
\quad \ 1       \quad\quad\quad\quad\quad\quad\quad \quad  \ \ \ \ \ \ \ \ \ \ \ \ \ \ \ \ \ \ \mbox{if}      \, \  |x|>1.
\end{cases}\end{equation*}
Then $\overline{x}=0$ is a tilt-stable local minimizer, and $f$  is  subdifferentially continuous but not prox-regular at $\overline{x}=0$ for $\overline{v}=0;$ see~\cite{DL} for further detail.  Here we check the prox-regular property of $f$ at $\ox$ for $\ov$ via definition and direct computation. Indeed, for each $\gamma>0$ sufficiently small,  we have
 $$M_\gg(v):={\rm argmin}\big\{f(x)-\la v,x\ra\big|\;x\in\B_\gg(0)\big\} =\{0\}\quad \mbox{for all}\  v\ \mbox{near} \ 0.$$
So, $\overline{x}=0$ is a tilt-stable local minimizer of $f$ with arbitrary  modulus $\kappa >0.$
Moreover, by simple calculations,  we get
$$\partial f(x) = \begin{cases}  \big[0,\frac{n+1}{n}\big]\ \ \ \ \ \ \ \mbox{if}   \  \    \  x=\frac{1}{n+1},\\
\big[-\frac{n+1}{n},0\big]\ \ \ \, \mbox{if}   \ \    \  x=-\frac{1}{n+1},\\
\big\{\frac{n+1}{n}\big\} \quad \ \ \ \ \ \ \mbox{if}   \ \ \  \  \frac{1}{n+1}<x<\frac{n+2}{(n+1)^2},\\
\big\{-\frac{n+1}{n}\big\} \quad \ \  \mbox{if}   \ \ \  \  -\frac{n+2}{(n+1)^2}<x<-\frac{1}{n+1},\\
 \big\{0,\frac{n+1}{n}\big\} \ \ \ \ \ \ \mbox{if}    \ \ \  \  x=\frac{n+2}{(n+1)^2},\\
 \big\{0, -\frac{n+1}{n}\big\} \ \ \, \ \mbox{if}    \ \ \  \  x=-\frac{n+2}{(n+1)^2},\\
 \{0\}\ \ \ \ \ \ \ \ \ \  \ \ \ \mbox{if}    \ \  \ \  \frac{n+2}{(n+1)^2}<|x|<\frac{1}{n},\\
 [-1,1] \ \ \ \ \ \ \ \ \  \mbox{if}    \ \  \ \  x=0
\end{cases}$$
with  $n\in \mathbb N^*.$
Let  $r$ be an arbitrary positive number.  For  $u_n:=\frac{1}{n+1}$,  $x_n:=\frac{1}{n+2}$ and $v_n=0\in\partial f(u_n),$ it holds that
$$\begin{array}{rl} &f(x_n)-f(u_n)-\la v_n,  x_n-u_n\ra + \frac{r}{2}|x_n-u_n|^2\\
&=\frac{1}{(n+1)(n+2)}\Big[\frac{r}{2(n+1)(n+2)}-1\Big]< 0 \,\  \mbox{for  all} \  n \  \mbox{large  enough},\end{array}$$
or equivalently,
$$f(x_n)<f(u_n)+\la v_n, x_n-u_n\ra - \frac{r}{2}|x_n-u_n|^2\,\  {\rm for \ all} \  n \ {\rm large \ enough.}$$
Therefore,  $f$ is not prox-regular at $\overline{x}=0$ for $\overline{v}=0.$
Next, we will show that \eqref{3.2} is invalid  for each  $\eta, \kk>0.$ To see this, note that for each $(u,v)\in~{\rm gph}\partial f$ with $u \in (-1,1),$ we have
$$T_{{\rm gph}\partial f}(u,v)=\begin{cases} \  \R\times \{0\}\  \ \ \ \ \ \ \ \ \ \ \ \ \ \ \ \ \ \ \ \ \ \ \ \ \mbox{if} \  \ \ \frac{1}{n+1}<|u|<\frac{n+2}{(n+1)^{2}}\ \mbox{or} \ \frac{n+2}{(n+1)^2}<|u|<\frac{1}{n},\\
\begin{array}{rl}\R_+\times \{0\}\  \ \ \ \ \ \ \ \ \ \ \ \ \ \ \ \ \ \ \ \ \ \ \mbox{if}  &u=\frac{n+2}{(n+1)^2},\, v=0 \ \mbox{or} \ u=-\frac{n+2}{(n+1)^2},\, v=-\frac{n+1}{n}, \end{array} \\
\begin{array}{rl}\R_-\times \{0\} \ \ \ \ \ \ \ \ \ \ \ \ \ \ \ \ \ \ \ \ \ \ \ \mbox{if}  &u=\frac{n+2}{(n+1)^2},\, v=\frac{n+1}{n} \ {\rm or} \ u=-\frac{n+2}{(n+1)^2},\, v=0,\end{array}\\
\begin{array}{rl}  (\R_+\times\{0\})\cup (\{0\}\times \R_-) \quad     \mbox{if} \ u=\frac{1}{n+1},\,  v=\frac{n+1}{n}\ \mbox{or}\   u=-\frac{1}{n+1},\, v=0,\end{array} \\
\begin{array}{rl}  (\{0\}\times \R_+)\cup (\R_-\times\{0\}) \ \ \      \mbox{if} \  u=\frac{1}{n+1},\, v=0\ \mbox{or}\ u=-\frac{1}{n+1},\, v=-\frac{n+1}{n},\end{array} \\
\begin{array}{rl} \{0\}\times \R   \ \ \ \ \ \ \ \ \ \ \ \ \ \ \ \ \ \ \ \ \ \ \ \ \ \mbox{if}  &u=\frac{1}{n+1},\, 0<v<\frac{n+1}{n}\\
&\mbox{or} \ u=-\frac{1}{n+1},\, -\frac{n+1}{n}<v<0,
\end{array}\\
\begin{array}{rl} \R_+\times \R   \ \ \ \ \ \ \ \ \ \ \ \ \ \ \ \ \ \ \ \ \ \ \ \ \  \,\mbox{if}\ u=0,\, 0<v<1, \end{array} \\
\begin{array}{rl} \R_-\times \R   \ \ \ \ \ \ \ \ \ \ \ \ \ \ \ \ \ \ \ \ \ \ \ \ \  \,\mbox{if}\ u=0,\, -1<v<0, \end{array} \\
\begin{array}{rl} (\R_+\times \R_+)\cup (\R_-\times\R_-) \ \ \   \,\mbox{if}\ u=0,\, v=0.\end{array}
\end{cases}$$
It follows that
$$D\partial f(u|v)(w)=\begin{cases}\begin{array}{rl}\{0\} \ \  \mbox{if}  &\frac{1}{n+1}<|u|<\frac{n+2}{(n+1)^{2}}\ \mbox{or} \ \frac{n+2}{(n+1)^2}<|u|<\frac{1}{n}\ {\rm and} \ w\in\R,\\
&\mbox{or} \ u=\frac{n+2}{(n+1)^2},\, v=0 \ \mbox{or} \ u=-\frac{n+2}{(n+1)^2},\, v=-\frac{n+1}{n}\  \mbox{and} \ w\geq0,\\
&\mbox{or}\  u=\frac{n+2}{(n+1)^2},\, v=\frac{n+1}{n} \ {\rm or} \ u=-\frac{n+2}{(n+1)^2},\, v=0 \ \mbox{and} \  w\leq0,\\
&\mbox{or}\  u=\frac{1}{n+1},\,  v=\frac{n+1}{n}\ \mbox{or}\   u=-\frac{1}{n+1},\, v=0 \ \mbox{and} \  w>0,\\
&\mbox{or}\  u=\frac{1}{n+1},\, v=0\ \mbox{or}\ u=-\frac{1}{n+1},\, v=-\frac{n+1}{n} \ \mbox{and} \  w<0,
\end{array}\\
\begin{array}{rl}\R_- \ \ \mbox{if} &u=\frac{1}{n+1},\,  v=\frac{n+1}{n}\ \mbox{or}\   u=-\frac{1}{n+1},\, v=0\ {\rm and} \ w=0,\\
&\mbox{or} \ u=0,\, v=0\  \mbox{and} \ w<0,\end{array}\\
\begin{array}{rl}\R_+ \ \  \mbox{if}  &u=\frac{1}{n+1},\, v=0\ \mbox{or}\ u=-\frac{1}{n+1},\, v=-\frac{n+1}{n}\ {\rm and} \ w=0,\\
&\mbox{or} \ u=0,\, v=0\  \mbox{and} \ w>0,\end{array}\\
\begin{array}{rl}\R \ \ \ \,     \mbox{if}  &u=\frac{1}{n+1},\, 0<v<\frac{n+1}{n}\ \mbox{or} \ u=-\frac{1}{n+1},\, -\frac{n+1}{n}<v<0, \ {\rm and} \ w=0,\\
&\mbox{or} \ u=0,\, 0<v<1\   \mbox{and} \ w\geq 0,\\
&\mbox{or}\  u=0,\, -1<v<0 \   \mbox{and} \ w\leq 0,\\
&\mbox{or}\  u=0,\, v=0 \   \mbox{and} \ w=0,
\end{array}\\
\ \ \emptyset \ \  \  \ {\rm otherwise.}
\end{cases}$$
For $(u_n,v_n):=(\dfrac{1}{n+1},0)\in \gph\partial f\cap\B_\eta(\ox,0)$ as $n\in\mathbb{N^\ast}, n>\dfrac{1}{\eta}-1$ , $w:=-1$ and $z:=0~\in ~D\partial f(u_n|v_n)(w),$ we have $\la z,w\ra=0<\frac{1}{\kk}\|w\|^2.$ Thus for any $\kk>0,$ the assertion (ii) in Theorem~\ref{thm1} is invalid.  This means that, in the absence of the prox-regularity,  tilt-stability of $f$ at  $\bar x$  does not guarantee the validity of  the second-order condition \eqref{3.2}.
}\end{Example}

\begin{Example}{\rm {\bf (Implication (ii) $\Rightarrow$ (i) fails in the absence of the prox-regularity).}\label{exam3.5}
Let  $f: \R^2 \rightarrow \overline{\R}$  be the function taken from \cite[Remark~4.8]{DMN}:
\begin{equation*}
f(x):= x_1^2+x_2^2+\delta_\Omega(x_1,x_2),
\end{equation*}
where $\Omega:=\{(x_1,x_2)\in\R^2\,| \, x_1x_2=0\}$ and $x=(x_1,x_2).$  We next show  that $f$ is not prox-regular at $\ox=0$ for $\overline{v}=0\in \partial f(0)$ and $\ox$ is not a tilt-stable local minimizer, while   the assertion {\bf (ii)} in Theorem~\ref{thm1} holds. To see this, we first note that, for each  $x=(x_1,x_2)\in\R^2,$  the limiting subdifferential of $f$ at $x$ is computed by
 $$\partial f(x_1,x_2)=(2x_1,2x_2)+N_\Omega(x_1,x_2)= \begin{cases}\{2x_1\}\times\R\ \ \  \ \mbox{if} \ \ \ x_1\not=0, x_2 =0,\\
\R\times\{2x_2\}\ \ \  \ \mbox{if} \ \ \ x_1=0, x_2\not =0,\\
\ \ \  \  \Omega  \ \ \ \ \ \  \ \ \ \ \ \mbox{if} \ \ \ x_1=0, x_2 =0,\\
\ \ \ \ \emptyset  \ \ \ \ \ \ \ \ \ \ \ \ \mbox{if} \ \ \ x_1x_2\not =0.
\end{cases}
$$
This implies that
$${\rm gph}\, \partial f =\big\{(x_1,0,2x_1,a)\, |\, x_1,a\in \R\big\} \cup \big\{(0,x_2,b,2x_2)|\; x_2,b\in \R\big\}.$$
For a fixed  $r>0,$   choosing  $u_n=\Big(\dfrac{1}{n},0\Big),$  $x_n=\Big(0,\dfrac{1}{n}\Big),$ and $v_n=\Big(\dfrac{2}{n},1\Big)\in\partial f(u_n),$ we have
$$f(x_n)-f(u_n)-\la v_n,x_n-u_n\ra+\frac{r}{2} \|x_n-u_n\|^2=(2+r-n)\frac{1}{n^2}<0 \ \mbox{as}\ n \ \mbox{large   enough.}$$
This says  $f$ is not prox-regular at $\ox=0$ for $\overline{v}=0.$  Moreover, by direct computation, for each $\gamma>0$ and $v=(v_1,v_2)$ with  $|v_1|,|v_2|<2\gamma,$ we obtain that
$$M_\gamma(v)=\begin{cases}\big\{\big(\frac{v_1}{2},0\big)\big\}\ \ \ \ \ \ \ \ \ \ \ \ \ \ \ \mbox{if} \ \ \ \ |v_1|>|v_2|,\\
\big\{\big(0,\frac{v_2}{2}\big)\big\}\ \ \ \ \ \ \ \ \ \ \ \ \ \ \ \mbox{if} \ \ \ \ |v_1|<|v_2|,\\
\big\{\big(\frac{v_1}{2},0\big),\ \big(0,\frac{v_2}{2}\big)\big\}\ \ \  \ \mbox{if} \ \ \ \ |v_1|=|v_2|.
\end{cases}$$
Since  $M_\gamma(v)$ is not single-valued,  $\ox$ is not a tilt-stable local minimizer.
Let  $(u,v)\in {\rm  gph}\partial f$ with  $u=(u_1,u_2)$ and $v=(v_1,v_2).$ For each  $w=(w_1,w_2)$ and $z=(z_1,z_2),$  by a simple calculation, we get $$T_{{\rm gph}\partial f}(u,v)=\begin{cases} \begin{array}{rl}\big\{(w,z)\in\R^2\times \R^2\,| \,w_2=0, \, z_1=2w_1\big\}  \ \ \mbox{if}  &u_1\not=0, u_2=0, \\
& \mbox{or}\  u_1= u_2= v_1=0,\, v_2\not=0,
\end{array}\\
\begin{array}{rl}\big\{(w,z)\in\R^2\times \R^2\, |\, w_1=0, \, z_2=2w_2\big\} \ \ \mbox{if} &u_1=0, u_2\not=0,\\
& \mbox{or}\  u_1= u_2= v_2=0,\, v_1\not=0,
\end{array}\\
\begin{array}{rl}&\big\{(w,z)\in\R^2\times \R^2\, |\, w_2=0, \, z_1=2w_1\big\}\\
&\cup \big\{(w,z)\in\R^2\times \R^2\, |\, w_1=0, \, z_2=2w_2\big\} \ \ \mbox{if} \ \ u_1= u_2= v_1=v_2=0.
\end{array}
\end{cases}$$
It follows that
$$D\partial f(u|v)(w)=\begin{cases} \begin{array}{rl} \{2w_1\}\times \R \ \ \ \ \ \ \mbox{if} &u_1\not=0, u_2=0, \\
& \mbox{or}\  u_1= u_2= v_1=0,\, v_2\not=0, \ {\rm and} \ w_2=0,
\end{array}\\
\begin{array}{rl} \R\times  \{2w_2\} \ \ \ \ \ \ \mbox{if}  &u_1=0, u_2\not=0,\\
& \mbox{or}\  u_1= u_2= v_2=0,\, v_1\not=0, \ {\rm and} \ w_1=0, \end{array}\\
\begin{array}{rl}\{2w_1\}\times \R \ \ \ \ \ \ \mbox{if}  &u_1= u_2= v_1=v_2=0 \ {\rm and} \ w_2=0,\, w_1\not=0, \end{array}\\
\begin{array}{rl} \R\times  \{2w_2\} \ \ \ \ \ \ \mbox{if}   &u_1= u_2= v_1=v_2=0 \ {\rm and} \ w_1=0,\, w_2\not=0, \end{array}\\
\begin{array}{rl} (\R\times  \{0\})\cup (\{0\}\times \R) \ \  \mbox{if}  &u_1= u_2= v_1=v_2= w_1=w_2=0, \end{array}\\
\ \ \emptyset \ \ \ \ \ \ \ \ \ \ \ \ \ \ \ \ \ \ \ \ \  \ \ \ \ \ \ \mbox{otherwise.}
\end{cases}$$
Hence, for  any   $(u,v)\in {\rm gph}\partial  f$ and $z \in  D\partial f(u|v)(w)$ with $w\in \R,$ we have
$ \langle z,w \rangle =2\|w\|^{2},$
which ensures \eqref{3.2}  for any  $\eta > 0$ and $\kk =\dfrac{1}{2}.$ We conclude that the validity of \eqref{3.2} does not imply tilt-stability of $f$ at $\ox$ without the prox-regularity assumption on $f$ at $\ox$ for $\ov$.

}\end{Example}

\section{Tilt Stability in  Nonlinear Programming}\label{TSNP}
\setcounter{equation}{0}

 In this section, using  our subgradient graphical derivative characterization of tilt-stability along with  the recent formulas for the graphical derivative of normal cone mappings from \cite{CH16, GY16} and some techniques from \cite{MG15}, we  establish  new results on tilt stability for nonlinear programming problems under the metric subregular constraint qualification.

  Consider the nonlinear programming problem:
  \begin{equation}\label{tiltMP} \begin{array}{rl} &\mbox{minimize} \quad g(x)\quad \mbox{subject to} \quad q_i(x)\leq 0,\ i=1,2,...,m,
 \end{array}
 \end{equation}
where $g:\R^n\rightarrow\R$ and $q_i:\R^n\rightarrow\R$  are twice continuously differentiable functions.

 Let $q: \R^n\rightarrow\R^m$  be the mapping defined  by $q(x):=\big(q_1(x),q_2(x),...,q_m(x)\big)$ for $x\in \R^n$ and $\Gamma:=\{x\,  |\,  q(x)\in\R^m_-\}$ be  the feasible set. Problem \eqref{tiltMP} could be written as a unconstrained optimization problem:
 \begin{equation*}\label{NP2} \mbox{minimize} \quad f(x):=g(x)+\delta_\Gamma(x),
 \end{equation*}
 where $\delta_\Gamma(x)$ is the indicator function to $\Gamma$, which equals to $0$ when $x\in \Gamma$ and $\infty$ otherwise. We say the point $\ox\in \Gamma$ is  a  {\it tilt stable local minimizer} of   Problem \eqref{tiltMP} with
modulus $\kk>0$  if there exists $\gg>0$ such that the argmin solution mapping
\[
\tilde M_\gg(v):=\disp{\rm argmin}\,\Big\{g(x)-\la v,x\ra|\; q(x)\in \R^m_+, x\in \overline\B_\gg(\ox)\Big\}
\]
is single-valued and Lipschitz continuous with constant $\kk>0$ on some neighborhood of $0\in \R^n$ with $\tilde M_\gg(0)=\ox$. Thus, $\ox$ is a  {\it tilt stable local minimizer} of  Problem \eqref{tiltMP} if and only if it  is a  {\it tilt stable local minimizer} of the function $f$ defined above. The number  ${\rm tilt}(g,q,\ox):={\rm tilt}(f,\ox)$ is the {\it exact modulus of tilt stability} of  \eqref{tiltMP} at  $\ox$.

 Following the sum rule \cite[Proposition~1.107]{M1}, the limiting subdifferential of $f$ at $x\in \Gamma$ is computed by
\begin{equation}\label{sumre}
\partial f(x)= \partial (g+\delta_\Gamma)(x)=\nabla g(x)+N_\Gamma(x):=\Psi(x)
\end{equation}
with $\Psi:\R^n\tto\R^n$.

Let us now recall some well-known constraint qualification in nonlinear programming.  The Mangasarian-Fromovitz constraint qualification (MFCQ) is said to hold  at the point $\ox\in \Gamma$ if there exists a vector $d\in \R^n$ such that
\[
\la \nabla q_i(\ox),d\ra<0\quad\mbox{for all}\quad  i\in I(\ox),
\]
where $I(\ox) :=\big\{i\in \{1, \ldots,m\}\, |\,  q_i(\ox)=0\big\}$ is the active index set at $\ox\in \Gamma$.  Furthermore,  the {\em constant rank constraint qualification} (CRCQ) is said to  hold at $\ox$ if there is a neighborhood $\mathcal{U}$ of $\ox$ such that the gradient system $\{\nabla q_i(x)|\;i\in J\}$ has the same rank in $\mathcal{U}$ for any index $J\subset I(\ox)$. It is well-recognized that MFCQ and CRCQ are independent in the sense that one cannot imply the other. Obviously,  CRCQ is weaker than  the {\em linear independence constraint qualification} (LICQ), which means all vectors $\nabla q_i(\ox)$, $i\in I(\ox)$ are linearly independent. Moreover, MFCQ is known to be stronger than the below {\em metric subregularity constraint qualification}; see, e.g.,  \cite{MG15}.

\begin{Definition}\label{Def41}{\rm

 ${\bf (i)}$ One says the {\it metric subregularity constraint qualification} (MSCQ) holds at $\bar x\in\Gamma$ if  the set-valued mapping $Q_q(x):=q(x)-\R^m_-$ is metrically subregular at $\bar x$ for  $0$, i.e., there exists a neighborhood $U$ of $\ox$ and a constant $\kk>0$ such that
 \begin{equation}\label{subre}
 d(x;\Gamma)=d(x;Q_q^{-1}(0))\le \kk\, d(0; Q_q(x))=\kk\, d(q(x);\R^m_-)\quad \mbox{for all}\  x\in U.
 \end{equation}
The infimum of all $\kk$ for which this inequality \eqref{subre}  holds is called  the {\it modulus of metric subregularity} of $Q_q$ at $\bar x$ for $0$ and is denoted by $ {\rm subreg}Q_q(\ox|0)$.

 ${\bf (ii)}$ The feasible set $\Gamma$ is said to have the {\it bounded extreme point property} (BEPP)  at $\ox\in \Gamma$ if there exist real numbers $\kk > 0$ and $r>0$  such that
  \begin{equation*}
 \mathcal{E}(x,x^\ast) \subset \kk\|x^\ast\|\B \quad \mbox{for  all}  \ x \in \Gamma \cap \B_r(\bar x)  \ \mbox{and} \ x^\ast \in \R^n,
 \end{equation*}
where  $\mathcal{E}(x,x^\ast)$ denotes  the set of extremal points of $\Lambda(x,x^\ast).$}
\end{Definition}

   The   metric subregularity constraint qualification  (MSCQ)  is a mild condition, which is  equivalent to  the existence of a local error bound~\cite{HO05}. It   is weaker than  most  known constraint qualifications, such as LICQ, MFCQ, CRCQ, the  pseudonormality and the quasinormality \cite{MT11}, the constant
positive linear dependence (CPLD)  \cite{QW00}, the relaxed CPLD \cite{AHSS11, GZL14}, the relaxed CRCQ \cite{MS11}, the relaxed MFCQ/the constant rank of the subspace component \cite{AHSS11b,  GZL14, KMO14}. Furthermore,   if MSCQ holds at $\bar x\in\Gamma$, then it  holds  at every $x\in \Gamma$ near $\bar x$. Note that MSCQ does not imply BEPP (see Example~\ref{vd11}), while the latter holds  under MFCQ or CRCQ or the  second-order sufficient condition for metric subregularity~\cite{MG15}.


By \cite[Theorem 31(b)]{CT10},  \cite[Proposition 3.4]{IO08} and  \cite[Corollary 1.15]{M1} and  one has the following result:
\begin{Lemma}\label{InvN} Let  $q:\R^n\rightarrow \R^m$ be  a  twice continuously differentiable mapping.
Suppose  MSCQ  holds at $\bar x\in \Gamma.$ Then,  there exists $\delta>0$ such that
\begin{equation}\label{InvNeq}N_\Gamma(x)=\widehat{N}_\Gamma(x)=\nabla q(x)^TN_{\R^m_-}( q(x))\ \, \mbox{for all}\ x\in\Gamma\cap\B_\delta(\bar x).\end{equation}
Moreover, $\delta_\Gamma$ is prox-regular and subdifferentially continuous at $\bar x$ for each $\bar x^*\in \partial\delta_\Gamma(\bar x).$
\end{Lemma}

Following Lemma~\ref{InvN}, under MSCQ, the normal cone to $\Gamma$ at $x$ could be presented by the following formula
\begin{equation}\label{noco}
N_\Gamma(x)=\Big\{\nabla q(x)^T\lambda\, |\, \lambda \in N_{\R^m_+}(q(x))\Big\}=\Big\{\nabla q(x)^T\lambda\, |\, \lambda\in \R^m_+, \lm_i=0 \mbox{ for } i\notin I(x)\Big\}.
\end{equation}
 For $x^*\in N_\Gamma(x)$, we denote
\[
 \Lambda(x,x^\ast):=\Big\{\lambda \in \R^m_+ \, |\, \nabla q( x)^T\lambda= x^*,\, \lm_i=0 \mbox{ for } i\notin I(x)\Big\}
\]
by the set of  Karush-Kuhn-Tucker (KKT) multipliers  corresponding to $(x,x^*)$. Furthermore, the critical cone to $\Gamma$ at $x$ for $x^*\in N_\Gamma(x)=[T_\Gamma(x)]^\circ$  is defined by
\[
K(x,x^*):=T_\Gamma(x)\cap \{ x^*\}^\bot.
\]
Let  $I^+(\lambda) :=\{i=1, \ldots,m\, |\, \lambda_i>0\}$ for $\lambda\in\R^m_+$. We note that if $\lm\in \Lm(x,x^*)$ then $I^+(\lambda)$ is a subset of the active set at $x$, i.e., $I^+(\lambda)\subset I(x)$.

  By Lemma~\ref{InvN},  if  MSCQ holds  at $\bar x,$ then the set  $\Lm(\ox,\ox^*)$ is a nonempty polyhedral convex set for any $\ox^*\in N_\Gamma(\ox)$. In this case, for each $v\in \R^n,$  the  problem
\[
{\rm LP}(v)\quad \quad\begin{array}{rl}  & \mbox{minimize} \quad  -v^T\nabla^2\big(\lambda^Tq\big)(\bar x)v\\
& \mbox{subject to}\quad \ \,   \lambda\in \Lm(\ox,\ox^*)\end{array}
\]
is a linear programming.
The optimal solution set of ${\rm LP}(v)$ will be denoted by $\Lambda(\bar x, \bar x^*;  v)$.

For problem \eqref{tiltMP}, its associated  KKT function is $L: \R^n\times \R^m\rightarrow \R$ defined by  $L(x,\lambda):=g(x)+\lambda^Tq(x)$  for each $x\in\R^n$ and $\lambda\in\R^m.$
Under MSCQ at a local minimizer $\ox$  to Problem \eqref{tiltMP}, it follows from \eqref{sumr} and \eqref{InvNeq} that there  exists  $\lambda\in \R^m_+$ such that  $\bar x$ is a solution of the KKT  system
\begin{equation}\label{n6.3}
\begin{cases} \nabla_xL(x,\lambda)=0,\\
 \lambda_iq_i(x)=0,\ i=1,2,...,m.\end{cases}
\end{equation}
When a  feasible point $x\in \Gamma$ satisfies~\eqref{n6.3} for some  KKT multiplier  $\lambda\in \R^m_+$, we call it a {\it stationary point} of  \eqref{tiltMP}.

In this paper, we introduce a new second-order sufficient condition, which is motivated from the so-called {\it uniform second-order sufficient condition} (USOSC) introduced by Mordukhovich and Nghia \cite{MN}.

 \begin{Definition} {\rm {\bf (Relaxed uniform second-order sufficient condition).} We say that the {\it relaxed  uniform second-order sufficient condition} (RUSOSC) holds at $\ox\in \Gamma$ with modulus $\ell>0$ if there  exists $\eta>0$ such that
\begin{equation}\label{usosc} \la  \nabla^2_{xx}L(x,\lambda)w,w\ra \geq \ell \|w\|^2\end{equation}
whenever $(x,v)\in {\rm gph}\Psi\cap \B_\eta(\ox, 0)$ with $\Psi$ defined in \eqref{sumre} and  $\lambda\in \Lambda\big(x,v-\nabla g(x); w\big)$ with $w\in \R^n$ satisfying
 \begin{equation}\label{Ron}
 \la \nabla q_i(x), w\ra=0\, \mbox{ for }\, i\in I^+(\lambda)\, \mbox{ and }\, \la \nabla q_i(x), w\ra \geq 0\, \mbox{ for }\, i\in I(x)\backslash I^+(\lambda).
 \end{equation}
}\end{Definition}
\begin{Remark}{\rm  The USOSC aforementioned in \cite{MN} is defined similarly, except for replacing   $\Lambda\big(x,v-\nabla g(x); w\big)$ by  $\Lambda\big(x,v-\nabla g(x)\big),$ which does not depend on $w$ satisfying \eqref{Ron}.  It is clear that USOSC implies  RUSOSC. The converse implication is also valid under CRCQ; see our Corollary~\ref{hq4.4} below together with Theorem~\ref{fnl-thm}.
}\end{Remark}

We now arrive at the first result of this section, which gives us a fuzzy characterization of tilt stable local minimizers in terms of  RUSOSC and its modification  for nonlinear programming problems. It is also worth noting here that the modulus of metric regularity ${\rm subreg}\,Q_q(\bar x|0)$ used in this result and the following ones could be computed directly in terms of initial data whenever MSCQ holds at $\ox$; see, e.g., \cite[Corollary 3.4]{MG17}.

\begin{Theorem}{\rm {\bf (Fuzzy characterization of tilt-stability under MSCQ).}} \label{fnl-thm} Given a stationary point  $\bar x\in \Gamma$ and  real numbers $\kk, \gamma>0,$ suppose that MSCQ is fulfilled at $\bar x$ and $\gamma>{\rm subreg}\,Q_q(\bar x|0).$  Then, the following assertions are equivalent:\label{thm2}

${\bf (i)}$  The point $\ox$ is a tilt-stable local minimizer of Problem \eqref{tiltMP} with modulus $\kk.$

${\bf (ii)}$ The RUSOSC is satisfied at $\ox$ with modulus $\ell:=\kk^{-1}.$

 {\bf (iii)}  There  exists $\eta>0$ such that
\begin{equation*}\label{usosc0} \la  \nabla^2_{xx}L(x,\lambda)w,w\ra \geq \frac{1}{\kk} \|w\|^2\end{equation*}
whenever $(x,v)\in {\rm gph}\Psi\cap \B_\eta(\ox, 0)$  and  $\lambda\in \Lambda\big(x,v-\nabla g(x); w\big)\cap \gamma\|v-\nabla g(x)\| \B_{\R^m}$ with $w\in \R^n$ satisfying
 \begin{equation*}\label{Ron0}
 \la \nabla q_i(x), w\ra=0\, \mbox{ for }\, i\in I^+(\lambda)\, \mbox{ and }\, \la \nabla q_i(x), w\ra \geq 0\, \mbox{ for }\, i\in I(x)\backslash I^+(\lambda),
 \end{equation*}
where $\Psi$ is defined in \eqref{sumre}.
\end{Theorem}
\noindent{\bf Proof.} Let $\eta > 0$ be so small that MSCQ holds at each  $x \in \Gamma\cap \B_{2\eta}(\overline{x})$ with modulus $\gamma.$  Pick any $z\in D\partial f(x,v)(w)$ with $(x,v)\in{\rm gph} \partial f\cap \B_{\eta}(\overline{x},0)={\rm gph} \Psi\cap \B_{\eta}(\overline{x},0)$ by \eqref{sumre}. It follows from the sum rule for the graphical derivative in \cite[Proposition 4A.2]{DR14} that
\begin{eqnarray}\label{sumr}D\partial f(x,v)(w)= \nabla^{2}g(x)(w)+DN_{\Gamma}(x,v-\nabla g(x))(w).\end{eqnarray}
By the computation of $DN_{\Gamma}(x,v-\nabla g(x))(w)$ in  \cite[Theorem 3.5]{CH16} and \cite[Theorem 4]{GY16}, we have
\begin{eqnarray}\label{comp} \begin{array}{rl}
& DN_{\Gamma}(x,v-\nabla g(x))(w)\disp= \bigcup_{\lambda\in \Lambda(x, v-\nabla g(x); w)} \nabla^{2}(\lambda^{T}q)(x)w+\widehat{N}_{K(x,v-\nabla g(x)) }(w)\\
&\disp =\bigcup_{\lambda\in \Lambda(x, v-\nabla g(x); w)\cap\gamma\|v-\nabla g(x)\| \B_{\R^m}} \nabla^{2}(\lambda^{T}q)(x)w +\widehat{N}_{K(x,v-\nabla g(x)) }(w).
\end{array}
\end{eqnarray}
Note further from the validity of  MSCQ, Lemma~\ref{InvN},  and \eqref{noco} that
\begin{equation}\label{tage}
T_{\Gamma}(x)=\big[\Hat N_\Gamma(x)\big]^\circ=\big\{u\in \R^{n}\, |\, \langle \nabla q_{i}(x),u\rangle \leq 0, \,  i \in I(x)\big\}.
\end{equation}
Moreover, for each $\lambda\in \Lambda\big(x,v - \nabla g(x)\big),$ we have
\[v - \nabla g(x) = \nabla q(x)^T\lambda = \sum\limits _{i\in I(x)}\lambda_{i}\nabla q_{i}(x)=\sum\limits _{i\in I^+(\lambda)}\lambda_{i}\nabla q_{i}(x),\]
which implies that
\begin{eqnarray}\label{crit}
w \in K(x, v -\nabla g(x))=T_{\Gamma}(x)\cap(v - \nabla g(x))^\perp \  \mbox{iff} \   \begin{cases}  \langle\nabla q_{i}(x),w\rangle = 0,
 i\in I^+(\lambda )\\
\langle\nabla q_{i}(x),w\rangle \le 0,   i\in  I(x) \backslash {I^+ }(\lambda ).
\end{cases}\end{eqnarray}

Let us justify  [${\bf (i)} \Rightarrow {\bf (ii)}$].  Assume that $\ox$ is a tilt-stable local minimizer of Problem~\eqref{tiltMP}, i.e., it is a tilt-stable local minimizer of $f$. Since $\eta>0$ could be arbitrarily small, we may suppose that \eqref{3.2} is satisfied with this $\eta$  by Theorem~\ref{thm1}.  Pick any $(x,v)\in {\rm gph}\partial f \cap \B_\eta(\ox, 0)= {\rm gph}\Psi\cap \B_\eta(\ox, 0),$ $w\in\R^n,$ and  $\lambda\in \Lambda\big(x,v-\nabla g(x); w\big)$ with
 \begin{equation*}
 \la \nabla q_i(x), w\ra=0\, \mbox{ for }\, i\in I^+(\lambda)\, \mbox{ and }\, \la \nabla q_i(x), w\ra \geq 0\, \mbox{ for }\, i\in I(x)\backslash I^+(\lambda).
 \end{equation*}
It follows from \eqref{crit} that $ -w \in K\big(x, v - \nabla g(x)\big)$.
Moreover, note from \eqref{InvNeq} and \eqref{eqdualTN} that
\[
v-\nabla g(x)\in N_\Gamma(x)=\Hat N_\Gamma(x)=[T_\Gamma(x)]^\circ\subset [K(x,v-\nabla g(x))]^\circ.
\]
 We have
\begin{equation}\label{nunu}
 \langle v - \nabla g(x),w\rangle =0 \ \, \mbox{and}\ \, v-\nabla g(x)\in\widehat{N}_{ K\big(x, v - \nabla g(x)\big)}(-w).
 \end{equation}
This implies $z:=-\nabla^2_{x}L(x,\lm)w+v-\nabla g(x)\in D\partial f(x,v)(-w)$ by \eqref{comp}.   Moreover, by Lemma~\ref{InvN} again, $f$ is prox-regular and subdifferentially continuous at $\bar x$ for $0.$ Thanks to \eqref{3.2} we obtain that
 \begin{eqnarray*}\big\langle-\nabla^{2}_{xx}L(x,\lambda)w+v-\nabla g(x),-w\ra=\langle z,-w\rangle\geq \frac{1}{\kappa} \|w\|^{2},
\end{eqnarray*}
This together with \eqref{nunu} verifies \eqref{usosc} with $\ell:=\kappa^{-1}$ and thus ensures {\bf (ii)}.

Since the implication [${\bf (ii)} \Rightarrow {\bf (iii)}$] is obvious,  it remains to  justify  [${\bf (iii)} \Rightarrow {\bf (i)}$]. To end this,  suppose   ${\bf (ii)}$ holds. Take any $(x,v)\in {\rm gph}\partial f \cap \B_{\eta}(\overline{x},0)={\rm gph} \Psi\cap \B_{\eta}(\overline{x},0)$ and $z\in D\partial f(x,v)(w).$
By \eqref{sumr} and \eqref{comp},  there exists $\lambda\in \Lambda\big(x,v-\nabla g(x); w\big)\cap \gamma\|v-\nabla g(x)\| \B_{\R^m}$ such that
 \begin{equation}\label{zzz}
  z-\nabla^{2}_{xx}L(x,\lm)w\in\widehat{N}_{K(x,v-\nabla g(x))}(w) \ \, \mbox{and}\,  \ w\in K(x, v-\nabla g(x)).
\end{equation}
It follows from  \eqref{crit} that $-w$ satisfies \eqref{Ron}. Noting also that $\Lambda\big(x,v-\nabla g(x); w\big)=\Lambda\big(x,v-\nabla g(x); -w\big).$ Hence, we get from ${\bf (iii)}$ and \eqref{usosc}  that
$$\la \nabla^2_{xx}L(x,\lambda)w,w\ra= \la -\nabla^2_{xx}L(x,\lambda)w,-w\ra\geq \ell \|-w\|^2=\ell \|w\|^2.$$
Moreover, since $K(x,v-\nabla g(x))$ is a cone, it follows from \eqref{zzz} that $\la z-\nabla^{2}_{xx}L(x,\lm)w,w\ra=0$.
Combining this with the above inequality gives us that
\[
\la z,w\ra=\la \nabla^2_{xx}L(x,\lambda)w,w\ra \geq \ell \|w\|^2.
\]
  By Theorem~\ref{thm1} and Lemma~\ref{InvN}, the point $\ox$ is a tilt stable local minimizer of \eqref{tiltMP} with modulus $\kk:=\ell^{-1}.$ The proof is complete.
\endproof

Using \cite[Proposition 5.3]{MG15}, we easily see that Theorem~\ref{fnl-thm} recovers  \cite[Theorem 4.3]{MN} under the validity of MFCQ and CRCQ.  Furthermore,  we show next that MFCQ is indeed  a superfluous assumption in  \cite[Theorem 4.3]{MN}.
\begin{Corollary}\label{hq4.4} {\rm {\bf (Characterization of tilt-stability under CRCQ via USOSC).}} Let $\ox$ be a stationary point of  \eqref{tiltMP} at which CRCQ holds. Then, the following assertions are equivalent: 

{\bf (i)}  The point $\ox$ is a tilt stable local minimizer of \eqref{tiltMP} with modulus $\kk>0.$

{\bf (ii)} There  exists $\eta>0$ such that
\begin{equation*} \la  \nabla^2_{xx}L(x,\lambda)w,w\ra \geq \frac{1}{\kk} \|w\|^2\end{equation*}
whenever $(x,v)\in {\rm gph}\Psi\cap \B_\eta(\ox, 0),$\, $\lambda\in \Lambda\big(x,v-\nabla g(x)\big),$  $\la \nabla q_i(x), w\ra=0$ for $i\in I^+(\lambda)$ and $\la \nabla q_i(x), w\ra \geq 0$ for $i\in I(x)\backslash I^+(\lambda).$
\end{Corollary}
\noindent{\bf Proof.}  Since CRCQ holds at $\bar x,$ by \cite[Proposition 5.3]{MG15},  we have
\[
\Lambda\big(x,v-\nabla g(x); w\big)=\Lambda\big(x,v-\nabla g(x)\big),
\]
 for all $(x,v)\in {\rm gph}\Psi$ near $(\bar x,0)$ and $w\in K(x,v-\nabla g(x)\big).$
Moreover,  the validity of CRCQ implies the validity of  MSCQ (see, e.g. \cite{MS11}), and  $\Lambda\big(x,v-\nabla g(x); w\big)=\Lambda\big(x,v-\nabla g(x); -w\big).$
Therefore, the conclusion follows from the equivalence between {\bf (i)} and {\bf (ii)} in Theorem~\ref{fnl-thm}. \endproof

 The below example shows a situation where Corollary~\ref{hq4.4} can be applicable, while  \cite[Theorem 4.3]{MN} cannot.

\begin{Example}{\rm  \label{examNCR}  Consider the following optimization problem in $\R^2$:
\begin{equation}\label{eq-ex4.1} \begin{array}{rl} & \mbox{minimize}\quad   x_1^2 + x_2^2,\\
&\mbox{subject to} \quad    x_{1}- x_{2} \leq 0,\   x_2-x_1 \leq 0.
\end{array}
\end{equation}
Obviously \eqref{eq-ex4.1}  is a special case of  \eqref{tiltMP} with   $g(x)=x_1^2 + x_2^2 ,$  $q_1(x)= x_1-x_2,$  $q_2(x)= x_2-x_1$  for $x=(x_1,x_2)\in \R^2$.  Note that  $\la w, \nabla^2_{xx}L(x,\lambda)w\ra = 2(w_{1}^2 + w_{2}^{2}),$ and  CRCQ holds  at $\bar x := (0,  0)$, so  does MSCQ.  On the other hand, direct computation shows that
 \[
 \la w, \nabla^2_{xx}L(x,\lambda)w\ra \geq 2\|w\|^2
 \]
 for all  $ v\in \Psi(x),$ $w\in\R^2,$ and  $ \lambda\in \Lambda\big(x,v-\nabla g(x); w\big)$ satisfying $$ \la \nabla q_i(x), w\ra=0\ \mbox{for}\  i\in I^+(\lambda),\ \mbox{and}\  \la \nabla q_i( x), w\ra \geq 0\   \mbox{for}\  i \in I(x)\backslash I^+(\lambda).$$
By  Corollary~\ref{hq4.4}, $\bar x= (0,  0) $ is  a tilt stable local minimizer for \eqref{eq-ex4.1}. However,  MFCQ does not hold at $\bar x,$  \cite[Theorem 4.3]{MN} is not applicable in this example.
}\end{Example}

The following result, which is a special case of \cite[Theorem 5.3.2~(i)]{BA}, is useful for us to obtain pointbased sufficient condition for tilt stability later.

\begin{Lemma}\label{fnl-lm1} {\rm (see \cite[Theorem 5.3.2 (i)]{BA}).} Consider the  problem $LP(a,c)$:
\begin{equation*}\label{tiltMP1} \begin{array}{rl}  &\mbox{\rm minimize} \ \quad a^Tx \\
& \mbox{\rm subject to}\quad b_i^T x \leq c_i , \   i = 1, \dots, m,
\end{array}
\end{equation*}
where  $ b_i\in\R^n, i=1,\dots ,m,$ $b=(b_1,...,b_m) \in \R^{mn}$  are given and fixed   and $a\in\R^n, c=(c_1,\dots ,c_m)\in\R^m.$ Let $\Phi(a,c)$ denote the set of   optimal  solutions to problem $ P(a,c).$ Then, the graph of the mapping  $\Phi: \R^n\times \R^m\rightrightarrows \R^n,$   $(a,c)\mapsto \Phi(a,c),$ is closed.
\end{Lemma}

The following result provides a pointbased sufficient condition for tilt-stable minimizer. 

\begin{Theorem}{\rm {\bf (Pointbased sufficient condition for tilt-stability  under MSCQ).}}  \label{fnl-thm11} Given a stationary point $\bar x\in \Gamma$ and  real numbers $\kk, \gamma>0,$ suppose that MSCQ is fulfilled at $\bar x$ and $\gamma>{\rm subreg}\,Q_q(\bar x|0)$   and that the following second-order condition holds:
\begin{eqnarray}\begin{array}{ll}\label{6.6}
&\qquad\qquad\qquad\qquad\la \nabla^2_{xx}L(\bar x,\lambda)w,w\ra>\frac{1}{\kk}\|w\|^2\\
&\mbox{whenever}\ w\not=0\ \mbox{with}\ \la \nabla q_i(\bar x), w\ra=0,\ i\in I^+(\lambda), {\rm and} \ \lambda\in\Delta(\ox),
\end{array}
\end{eqnarray}
where $\Delta(\ox):=\bigcup\limits_{0\not=v\in K\big(\bar x,-\nabla g(\bar x)\big)}\Lambda\big(\bar x, -\nabla g(\bar x); v\big)\bigcap  \gamma\|\nabla g(\bar x)\| \B_{\R^m}.$\\
Then $\bar x$ is a tilt-stable local minimizer of \eqref{tiltMP} with modulus $\kk.$ Furthermore, we
have the  estimation:
\begin{equation} \label{6.8}
{\rm tilt}(g,q,\ox)\leq\sup\bigg\{\dfrac{\|w\|^2}{\la  \nabla^2_{xx}L(\bar x,\lambda)w,w\ra}| \ \lambda\in\Delta(\ox), \la \nabla q_i(\bar x), w\ra=0,\ i\in I^+(\lambda)\bigg\}<\infty
\end{equation}
with the convention that $0/0 := 0$ in
\eqref{6.8}.
\end{Theorem}
\noindent{\it Proof.} Suppose to contrary that all assumptions  of  Theorem~\ref{fnl-thm11} are satisfied, but  $\bar x$ is not a tilt-stable local minimizer of \eqref{tiltMP} with modulus $\kk.$  By the equivalence between {\bf (i)} and {\bf (iii)} in Theorem~\ref{fnl-thm},  there exist $(x^k,v^k)\rightarrow (\bar x, 0)$ with   $v^k\in \Psi(x^k)=\nabla g(x^k)+N_\Gamma(x^k)$ and
$\lambda^k\in \Lambda\big(x^k,v^k-\nabla g(x^k); w^k\big)$ for some $w^k\in \R^n$ satisfying  $\|\lambda^k\|\leq \gamma\|v^k-\nabla g(x^k)\|$,  \begin{equation}\label{Ron00}
 \la \nabla q_i(x^k), w^k\ra=0\, \mbox{ for }\, i\in I^+(\lambda^k)\, \mbox{ and }\, \la \nabla q_i(x^k), w^k\ra \geq 0\, \mbox{ for }\, i\in I(x^k)\backslash I^+(\lambda^k),
 \end{equation}
and  that
\begin{equation}\label{usosc00} \la  \nabla^2_{xx}L(x^k,\lambda^k)w^k,w^k\ra < \frac{1}{\kk} \|w^k\|^2.\end{equation}
By dividing all equalities and inequalities in \eqref{Ron00} by $\|w^k\|$ and both sides of \eqref{usosc00}  by $\|w^k\|^2$,  we may assume without loss of generality that  $\|w^k\|=1$ and $(w^k)$ converges to some $\bar w\in \R^n$ with $\|\bar w\|=1.$  Furthermore, since  $\|\lambda^k\|\leq \gamma\|v^k-\nabla g(x^k)\|$ for all $k\in\mathbb N$ and $(x^k,v^k)\rightarrow (\bar x, 0),$  using a subsequence if necessary, we may assume  $\lambda^k\rightarrow \bar \lambda\in \R^m$ with $\|\bar\lambda\|\leq  \gamma\|\nabla g(\bar x)\|.$
By passing  $k\rightarrow \infty$ in \eqref{usosc00}, we get
 \begin{equation}\label{6.74} \left\langle \nabla_{xx}^2L(\bar x,\bar \lambda)\bar w, \bar w\right\rangle\leq\frac{1}{\kk}\|\bar w\|^2.\end{equation}
Note that $\lambda^k\in N_{\R^m_-}(q(x^k))$, $\lm^k\to \bar\lm$, $q(x^k)\to q(\ox)$, we get that $\bar \lambda\in N_{\R^m_-}(q(\bar{x}))$. Moreover, we have
\[
\nabla q(\bar{x})^T\bar\lambda=\lim_{k\rightarrow\infty}\nabla q(x^k)^T\lambda^k=\lim_{k\rightarrow\infty}\big(v^k-\nabla g(x^k)\big)=-\nabla g(\bar{x}).
\]
Hence,  $\bar \lambda\in \Lambda\big(\bar x, -\nabla g(\bar x)\big).$
We next show that there exists $v\in K\big(\bar x,-\nabla g(\bar x)\big)\backslash\{0\}$ such that      $\bar \lambda\in \Lambda\big(\bar x, -\nabla g(\bar x); v\big).$  Consider  the following two cases:

{\it Case 1:} There exist infinitely many $k\in \mathbb N$ such that  $x^k \neq\bar{x}.$  Passing to a subsequence if necessary, we may assume that   $x^k \neq\bar{x}$ for every $k$ and  that  $(x^k-\bar{x})/\|x^k-\bar{x}\|\rightarrow v$ for some $v\in \R^n$ with $\|v\|=1.$  Note that $I^+(\bar{\lambda}) \subset I^+(\lambda^k)\subset I(x^k)\subset I(\ox)$ for large $k$, we get
\[\langle\nabla q_i(\bar{x}),v\rangle=\lim_{k\rightarrow\infty}\dfrac{q_i(x^k)-q_i(\bar{x})}{\|x^k-\bar{x}\|}\begin{cases}=0\ \ \mbox{if} \ \ \  i \in I^+(\bar{\lambda}),\\
\leq 0 \ \ \mbox{if} \ \ \ i \in I(\bar x)\setminus I^+(\bar{\lambda}). \end{cases}\]
Hence, $v\in T_\Gamma(\ox)=\{w|\; \nabla q_i(\ox)w\le 0, i\in I(\ox)\}$ by \eqref{tage}. Moreover, since  $-\nabla g(\ox)=\nabla q(\ox)^T\bar\lm$, we get from the above expression that
\[
\la v, \nabla g(\ox)\ra=\sum_{i=1}^m\lm_i\la v,\nabla q_i(\ox)\ra=\sum_{i\in I^+(\bar \lm)}\lm_i\la v,\nabla q_i(\ox)\ra=0.
\]
It follows that
 $v\in T_\Gamma(\ox)\cap\{-\nabla g(\ox)\}^\perp= K\big(\bar{x},-\nabla g(\bar{x})\big).$  Note further that   $\left\langle\bar{\lambda},q(x^k)\right\rangle = 0$ for all  $k\in \mathbb N$  sufficiently large due to the fact $I^+(\bar \lm)\subset I(x^k)$. Pick any $\lm\in \Lm(\ox,-\nabla g(\ox))$, we have   $\left\langle\lambda,q(x^k)\right\rangle\leq 0$ and   $\left\langle\lambda,\nabla q(\bar{x})\right\rangle=\left\langle\bar{\lambda},\nabla q(\bar{x})\right\rangle=0.$ Therefore,
\[\begin{array}{rl} 0&\geq2\limsup\limits_{k\rightarrow\infty}\dfrac{\left\langle\lambda-\bar{\lambda},q(x^k)\right\rangle}{\|x^k-\bar{x}\|^2}
=2\limsup\limits_{k\rightarrow\infty}\dfrac{\left\langle\lambda-\bar{\lambda},q(x^k)\right\rangle-\left\langle\lambda-\bar{\lambda},q(\bar{x})\right\rangle}{\|x^k-\bar{x}\|^2}\\
&=\limsup\limits_{k\rightarrow\infty}\dfrac{\left\langle x^k-\bar{x},\nabla^2\langle\lambda-\bar{\lambda},q\rangle(\bar{x})(x^k-\bar{x})\right\rangle}{\|x^k-\bar{x}\|^2}=\left\langle v,\nabla^2\langle\lambda-\bar{\lambda},q\rangle(\bar{x})v\right\rangle, \end{array}\]
which clearly implies that $\bar \lambda\in\Lambda\big(\bar{x},-\nabla g(\bar{x});v\big)$.

{\it Case 2:}  $x^k\neq\bar{x}$ for  only finitely many $k.$  Then we  may assume  without loss of generality  that $x^k=\bar{x}$ for all $k\in \mathbb{N}.$  From \eqref{Ron00} it follows that
$$w^k\in K\big(\bar x,v^k-\nabla g(\bar x)\big)=T_\Gamma(\bar{x})\cap \{v^k-\nabla g(\bar x)\}^\bot.$$
So, since  $\|w^k\|=1$ and $(w^k,v^k)\rightarrow (\bar w, 0),$  we have $\bar{w}\in K\big(\bar x,-\nabla g(\bar x)\big)\backslash \{0\}.$  On the other hand,  $\lambda^k\in \Lambda\big(\bar x,v^k-\nabla g(\bar x); w^k\big).$  By Lemma~\ref{fnl-lm1},  we have $\bar{\lambda}\in\Lambda(\bar{x}, -\nabla g(\bar{x}); \bar{w}).$
Moreover, since    $\la \nabla q_i(x^k), w^k\ra=0$ for all $i\in I^+(\lm^k),$  we see that  $\la \nabla q_i(\bar x),\bar w\ra=0$ for all $i\in I^+(\bar\lambda).$

Consequently, in  the both cases above, we get a contradiction by comparing  \eqref{6.74} with  \eqref{6.6}. So, $\bar x$ is a tilt-stable local minimizer of \eqref{tiltMP} with modulus~$\kk.$
Finally, to justify \eqref{6.8}, note the assumption \eqref{6.6} that
\[
\rho:=\sup\bigg\{\dfrac{\|w\|^2}{\la  \nabla^2_{xx}L(\bar x,\lambda)w,w\ra}| \ \lambda\in\Delta(\ox), \la \nabla q_i(\bar x), w\ra=0,\ i\in I^+(\lambda)\bigg\}<\infty.
\]
Take any $\kk>\rho$, we have \eqref{6.6}. This along with the above proof guarantees that  $\bar x$ is a tilt-stable local minimizer of \eqref{tiltMP} with modulus $\kk,$ showing that   ${\rm tilt}(g,q,\ox)\le \kk$. Since $\kk>\rho$ is chosen arbitrarily, we have ${\rm tilt}(g,q,\ox)\le \rho$, which clearly verifies   \eqref{6.8} holds. The proof is complete.
   \endproof

As a consequence of Theorem~\ref{fnl-thm11},  we establish the following result, which is also a corollary of \cite[Theorem 6.1]{MG15} by replacing MSCQ and BEPP there by the stronger constraint qualification that is  MFCQ. 

\begin{Corollary}  \label{c1}
Let  $\bar x\in \Gamma$  be a stationary point of Problem \ref{tiltMP} and   $\kk$ be  a positive number. Suppose that MFCQ is  satisfied at $\bar x$   and the following second-order condition holds:
\begin{eqnarray}\begin{array}{ll}\label{6.6-1}
&\qquad\qquad\qquad\qquad\la \nabla^2_{xx}L(\bar x,\lambda)w,w\ra>\frac{1}{\kk}\|w\|^2\\
&\mbox{whenever}\ w\not=0\ \mbox{with}\ \la \nabla q_i(\bar x), w\ra=0,\ i\in I^+(\lambda), {\rm and} \ \lambda\in\Delta_{\mathcal{E}}(\ox),
\end{array}
\end{eqnarray}
where $\Delta_{\mathcal{E}}(\bar x):=\bigcup\limits_{0\not=v\in K\big(\bar x,-\nabla g(\bar x)\big)}\Lambda_{\mathcal{E}}\big(\bar x, -\nabla g(\bar x); v\big)$ and  $\Lambda_{\mathcal{E}}\big(\bar x, -\nabla g(\bar x); v\big)$ is the set of extremal points of  $\Lambda\big(\bar x, -\nabla g(\bar x); v\big).$
Then $\bar x$ is a tilt-stable local minimizer of \eqref{tiltMP} with modulus $\kk.$ Furthermore, we
have the estimate
\begin{equation} \label{6.8-1}
{\rm tilt}(g,q,\ox)\leq\sup\bigg\{\dfrac{\|w\|^2}{\la  \nabla^2_{xx}L(\bar x,\lambda)w,w\ra}| \ \lambda\in\Delta_{\mathcal{E}}, \la \nabla q_i(\bar x), w\ra=0,\ i\in I^+(\lambda)\bigg\}
\end{equation}
 with the convention that $0/0 := 0$ in
\eqref{6.8-1}.
\end{Corollary}
\noindent{\it Proof.}  Let  any $\lm\in \Lambda\big(\bar x, -\nabla g(\bar x); v\big)$  for some $v\in K\big(\bar x,-\nabla g(\bar x)\big)$ and $w\not=0$ with  $\la\nabla q_i(\ox),w\ra=0, \ i\in I^+(\lambda).$ Since MFCQ holds at $\bar x,$   the set $\Lambda(\ox,-\nabla g(\bar x);v)$ is bounded and thus it is  a compact polyhedral set. Hence, we have   $$\lambda\in  {\rm conv}\big\{\lambda\, |\, \lambda\in\Lambda_{\mathcal{E}}\big(\bar x, -\nabla g(\bar x); v\big)\big\}.$$
Consequently,   $\lambda=\sum\limits_{i=1}^s\mu_i\lambda^i$ for some  $\lambda^i\in  \Lambda_{\mathcal{E}}\big(\bar x, -\nabla g(\bar x); v\big),$ $\mu_i>0,$ $i=1,...,s,$ $s\in \mathbb N^*,$ and $\sum\limits_{i=1}^s\mu_i=1.$
Since $I^+(\lambda^i)\subset I^+(\lambda)$ and  $\la\nabla q_i(\ox),w\ra=0,$ $i\in I^+(\lambda),$ we get
$$\la\nabla q_i(\ox),w\ra=0, \ i\in I^+(\lambda^i)\quad \mbox{for all}\  i=1,...,s.$$
So, taking into account that $\lambda^i\in  \Lambda_{\mathcal{E}}\big(\bar x, -\nabla g(\bar x); v\big),$ by \eqref{6.6-1}, we have
$$\la w,\nabla^2_{xx}L(\ox,\lambda^i)w\ra>\frac{1}{\kk}\|w\|^2\quad \mbox{for}\ i=1,...,s.$$
This implies that
\begin{equation*}\label{eqce1} \la w,\nabla^2_{xx}L(\ox,\lambda)w\ra=\sum\limits_{i=1}^s\mu_i\la w,\nabla^2_{xx}L(\ox,\lambda^i)w\ra>\sum\limits_{i=1}^s\mu_i\frac{1}{\kk}\|w\|^2=\frac{1}{\kk}\|w\|^2.
\end{equation*}
  So, by  Theorem~\ref{fnl-thm11}, $\bar x$ is a tilt-stable local minimizer of \eqref{tiltMP} with modulus $\kk,$ and thus, \eqref{6.8-1} follows. \endproof

We next establish another second-order sufficient condition for tilt-stable local minimizers by surpassing the appearance $\kk$ in \eqref{6.6}.

\begin{Theorem}\label{cofill} Given a stationary point $\bar x\in \Gamma$ and  a real number $\gamma>0,$ suppose that MSCQ is fulfilled at $\bar x$ and $\gamma>{\rm subreg}\,Q_q(\bar x|0)$   and that the following second-order condition holds:
\begin{eqnarray}\label{613}\begin{array}{ll}
 &\la w, \nabla^2_{xx}L(\bar x,\lambda)w\ra>0\quad \mbox{whenever}\ w\not=0\ \mbox{with}\ \la \nabla q_i(\bar x), w\ra=0,\ i\in I^+(\lambda), \\
& and \ \  \ \ \ \lambda \in \Delta(\ox):=\bigcup\limits_{0\not=v\in K\big(\bar x,-\nabla g(\bar x)\big)}\Lambda\big(\bar x, -\nabla g(\bar x); v\big)\bigcap  \gamma\|\nabla g(\bar x)\| \B_{\R^m}.
\end{array}
\end{eqnarray}
Then,   $\ox$ is a tilt-stable local minimizer for \eqref{tiltMP}.
\end{Theorem}
\noindent{\it Proof.} If  $\Delta(\ox)=\emptyset$ or $\bigcup\limits_{\lambda\in \Delta(\ox)}\big\{ w\ \langle \nabla q_i(\ox),w\rangle =0,\ i\in I^{+}(\lambda)\big\}=\{0\},$ then, by Theorem~\ref{fnl-thm11}, we get the conclusion. We now assume  $\Delta(\ox)\not =\emptyset$ and $\bigcup\limits_{\lambda\in \Delta(\ox)}\big\{ w\ \langle \nabla q_i(\ox),w\rangle =0,\ i\in I^{+}(\lambda)\big\}\not=\{0\}.$    First, we justify the compactness of  $\Delta(\ox).$   Since $\Delta(\ox)$ is bounded, it suffices to prove that $\Delta(\ox)$ is closed.  To do this,  take any $\{\lambda^k\}\subset\bigcup\limits_{0\not=v\in K\big(\bar x,-\nabla g(\bar x)\big)}\Lambda\big(\ox, -\nabla g(\ox); v\big)$ with  $\lambda^k\rightarrow \overline{\lambda}.$ Since $\Lambda\big(\ox, -\nabla g(\ox); tv\big)=\Lambda\big(\ox, -\nabla g(\ox); v\big)$ for all $t>0, $ $v\in K\big(\ox,-\nabla g(\ox)\big)\backslash\{0\},$ and $K\big(\ox,-\nabla g(\ox)\big)$ is a cone,   one can find  $v_k\in K\big(\ox,-\nabla g(\ox)\big)\backslash\{0\}$ with $\|v_k\|=1$ such that  $\lambda^k\in \Lambda\big(\ox, -\nabla g(\ox); v_k\big)$ for all $k$.  By passing to subsequences if necessary,  we may assume that $v_k\rightarrow \overline{v}$ with $\|\overline{v}\|=1.$ Applying Lemma~\ref{fnl-lm1} to the situation of problem ${\rm LP}(v)$,  we have $\overline{\lambda}\in \Lambda\big(\ox, -\nabla g(\ox); \overline{v}\big).$  It follows  that $\Delta(\ox)$ is closed and thus is compact.

 Next, we show  that  \eqref{6.6} holds for some $\kk>0$ when \eqref{613} is satisfied.   Indeed, for each $\lambda\in\Delta(\ox)$, let
$$V_\lambda:= \{w\in\R^n| \la w,\nabla q_i(\ox)\ra=0, i\in I^+(\lambda)\},$$
$$\delta_\lambda:=\inf\{w^T\nabla^2L(\bar x,\lambda)w\, |\, \|w\|=1,\, \la w,\nabla q_i(\ox)\ra=0, i\in I^+(\lambda)\}$$
and   $$\ell:= \inf\limits_{\lambda\in\Delta(\ox),\, V_\lambda\not=\{0\}}\delta_\lambda.$$ We  note that  \eqref{613} ensures $\delta_\lambda>0$  for all $\lambda\in\Delta(\ox)$ with $V_\lambda\not=\{0\}.$ From the definition of $\ell$ there exists $\{\lambda^k\}\subset\Delta(\ox)$ with $V_{\lambda^k}\not=\{0\}$ such that $\lim\limits_{k\rightarrow\infty}\delta_{\lambda^k}=\ell.$  Since $\Delta(\ox)$ is compact and the sets $I(\lambda)$ has finite elements in $\{1,2\ldots,m\}$, $\lm\in \Delta(\ox)$,  passing to a subsequence if necessary,  we may  assume that $\lambda^k\rightarrow\overline{\lambda}$ for some $ \overline{\lambda} \in \Delta(\ox),$
 and   $I^+(\lambda^k)=I^+(\lambda^1)\supset I^+(\bar \lambda)$  for all $k.$  Therefore,  $ V_{\lambda^k}=V_{\lambda^1}\subset V_{\overline{\lambda}}$ for all $k.$ This implies that  $V_{\overline{\lambda}}\not=\{0\}$ and  $$\nabla^2L(\ox,\lambda^k) \rightarrow \nabla^2L(\ox,\overline{\lambda})\,\ \mbox{as}\ \, k\rightarrow\infty.$$
 We see that
 $$\ell=\lim\limits_{k\rightarrow\infty}\delta_{\lambda^k}=\inf\{w^T\nabla^2L(\bar x,\bar\lambda)w\, |\, \|w\|=1,\, \la w,\nabla q_i(\ox)\ra=0, i\in I^+(\lambda^1)\}.$$
 So, taking into account that  $I^+(\lambda_1)\supset I^+(\bar \lambda),$ we have
 $$\begin{array}{rl}\ell&\geq \inf\{w^T\nabla^2L(\bar x,\bar\lambda)w\, |\, \|w\|=1,\,  \la w,\nabla q_i(\ox)\ra=0, i\in I^+(\bar \lambda)\}\\
 &=\delta_{\bar\lambda}>0.\end{array}$$
 Finally, for each $\lambda\in\Delta(\ox)$ and  $ 0\not=w\in V_\lambda,$  we have
$$\begin{array}{rl} \la w,\nabla^2_{xx}L(\ox,\lambda)w\ra &\geq \delta_\lambda\|w\|^2\\
&\geq\ell\|w\|^2.\end{array}$$
 This shows that \eqref{6.6} holds for any $\kappa>\ell^{-1}$. By Theorem~\ref{fnl-thm11}, $\ox$ is a tilt-stable local minimizer with modulus $\kappa.$   \endproof

Recall that the {\it strong second-order sufficient condition} (SSOSC) holds at $\ox$ if for all $\lambda\in\Lambda\big(\ox,-\nabla g(\ox)\big)$ we have
\begin{equation}\label{SSOSC}
\la w, \nabla^2_{xx}L(\bar x,\lambda)w\ra>0\quad \mbox{whenever}\ w\not=0\ \mbox{with}\ \la \nabla q_i(\bar x), w\ra=0,\ i\in I^+(\lambda).
\end{equation}
Under MFCQ and CRCQ, Mordukhovich and Outrata \cite[Theorem 3.5]{MO1} proved that  the tilt-stability is satisfied under SSOSC. In the following corollary  we also obtain this property but under weaker condition.

\begin{Corollary} \label{ssosc}{\rm {\bf (Tilt stability from SSOSC under MSCQ).}} Let  $\bar x$  be a stationary point of  \eqref{tiltMP} at which  MSCQ is valid.  Then, $\bar x$ is a tilt-stable local minimizer of  \eqref{tiltMP} provided SSOSC is satisfied  at $\bar x.$
  \end{Corollary}
\noindent{\it Proof.} The desired conclusion is straightforward from   Theorem~\ref{cofill}. \endproof

To complete this section, we provide a simple example, which is accessible  by our Theorem~\ref{fnl-thm11} and Theorem~\ref{cofill} to verify the tilt stability; however, all the results in  \cite{MG15} is not applicable, since BEPP is not valid in this example.  Moreover, we also clarify Corollary~\ref{ssosc} without using either MFCQ and CRCQ as in \cite[Theorem~3.5]{MO1} discussed above.

\begin{Example}\label{vd11}{\rm {\bf (Tilt stability under MSCQ without BEPP).}
Consider the following two-dimensional nonlinear programming problem:
\begin{equation}\label{vd1} \begin{array}{rl} &\mbox{minimize} \ \quad x_2^2+x_1x_2-x_1 \\
 & \mbox{subject to}\quad -x_1\leq 0,\ x_1\leq 0, \  x_1x_2^2\leq 0.
 \end{array}
 \end{equation}
 Let  $g(x): = x_2^2+x_1x_2-x_1,$ $q_1(x):=-x_1,$ $q_2(x)=x_1,$ $q_3(x)=x_1x_2^2,$    $q(x):=\big(q_1(x),q_2(x), q_3(x)\big),$ $x=(x_1,x_2)\in \R^2,$ and   $\bar x=(0,0)\in\Gamma.$  We have
 \[
 \Gamma:=\big\{x\in\R^2 \, |\, q(x)\in \R^3_-\big\}=\{0\}\times\R
 \]
  and $\bar x\in\Gamma.$ Furthermore, direct verification  shows  that
 \[
 \nabla q(x)=\begin{pmatrix}-1 & 1 & x_2^2 \\
 0 & 0 & 2x_1x_2
 \end{pmatrix}^T, \ \nabla q(\ox)=\begin{pmatrix}-1 & 1 & 0 \\
 0 & 0 & 0
 \end{pmatrix}^T,\  \nabla^2(\lambda^Tq)(\bar x)=\begin{pmatrix}0&0\\
 0&0 \end{pmatrix} \mbox{for}\ \lm\in \R^3, $$
  $$   \nabla g(\ox)=(-1,0), \ \nabla^2g(\ox)=\begin{pmatrix}0& 1\\
 1& 2
\end{pmatrix},\ K\big(\ox,-\nabla g(\ox)\big)=\{0\}\times\R,
\]
and
\[
 \Lambda\big(\ox,-\nabla g(\ox);v\big)=\Lambda\big(\ox,-\nabla g(\ox)\big) =\{\lambda\in\R^3_+\ |\ -\lambda_1+\lambda_2=1\}
 \]
  for each $v\in K\big(\ox,-\nabla g(\ox)\big),$ where $x=(x_1,x_2)\in\R^2$ and  $\lambda=(\lambda_1,\lambda_2,\lambda_3)\in\R^3.$ Obviously, we see that
  \[ d(x;\Gamma)=|x_1|\le \max\{-x_1,0\}+\max\{x_1,0\}+\max\{x_1x_2^2, 0\}=d\big(q(x); \R^3_-\big)
  \]
for all  $x=(x_1,x_2)\in \R^2.$ Hence  MSCQ is fulfilled  at $\ox$ and ${\rm subreg}\,Q_q(\bar x|0)=1,$ where $Q_q(x):=q(x)-\R^3_-.$
Let $\gamma>1$ and  $\Delta(\ox):=\bigcup\limits_{0\not=v\in K\big(\bar x,-\nabla g(\bar x)\big)}\Lambda\big(\bar x, -\nabla g(\bar x); v\big)\bigcap  \gamma\|\nabla g(\bar x)\| \B_{\R^3}.$
It is easy to see that
\[
\gamma>{\rm subreg}\,Q_q(\bar x|0)\quad \mbox{and}\quad  \Delta(\ox)=\{\lambda\in\R^3_+\ |\ -\lambda_1+\lambda_2=1, \ \|\lambda\|\leq \gamma\}.
\]
Take an arbitrary    $\lambda=(\lambda_1,\lambda_2,\lambda_3)\in\Delta(\ox).$ We note that   $\lambda_2>0.$ So, if  $\la \nabla q_i(\bar x), w\ra=0,$  $i\in I^+(\lambda),$ and  $w=(w_1,w_2)\not=0$, then  $w_1=0,$ $ w_2\not=0.$
Therefore,
\begin{equation}\label{zy}
\la w, \nabla^2_{xx}L(\bar x,\lambda)w\ra=2\|w\|^2\quad \mbox{for all}\ w\not=0,\  \la \nabla q_i(\bar x), w\ra=0,\ i\in I^+(\lambda),  \ \lambda\in\Delta(\ox),
\end{equation}
where $L(x,\lambda):=g(x)+\lambda^Tq(x)$ and $\kk>\frac{1}{2}.$  By Theorem~\ref{fnl-thm11}, $\bar x$ is a tilt-stable local minimizer of \eqref{vd1} with modulus $\kk.$

We next show that  BEPP does not hold at $\ox$ in this case.  Indeed, for each $i=1,2,...,$  letting  $x^i=(0,\frac{1}{i})\in \Gamma$ and  $ x^\ast:=-\nabla g(\bar x)=(1,0),$ we have
\[
\Lambda(x^i,x^\ast)=\big\{\lambda\in\R^3_+ \ |\ - \lambda_1+\lambda_2+\frac{1}{i^2}\lambda_3=1\big\}.
\]
 Thus  $\lambda^i:=(0,0,i^2)$ is a point in  $\mathcal{E}(x^i,x^\ast)$ with   $ \|\lambda^i\|\rightarrow\infty.$  This infer that  BEPP is not fulfilled at $\bar x.$ Therefore, there is no any result of~\cite{MG15} that can  apply to this example.

Finally, we observe that it is similar to \eqref{zy}  that SSOSC \eqref{SSOSC}  is satisfied at $\ox$, while either MFCQ or CRCQ fails in this example. This is an evidence of the advantage of our  Corollary~\ref{ssosc} in comparison to \cite[Theorem~3.5]{MO1}.

}\end{Example}


\section{Concluding Remarks}

In this paper we have introduced a new fuzzy characterization of tilt stability via the sugradient graphical derivative. This new approach allows us to  obtain some  second-order necessary and sufficient conditions for tilt stability in nonlinear programming,  which extend  and improve several recent  results in  \cite{MG15,MN,MO1} by weakening the involved assumptions.   Keeping in mind that,  in the current stage,  the commonly used dual approach has met severe difficulties in
handling tilt stability for non-polyhedral conic programs under weak conditions,  examining the new approach to tilt stability  for such problems   would be  a  topic of  great interest.
Another important topic of further research is  to expand our approach to {\em full stability} in the sense of Levy-Poliquin-Rockafellar \cite{LPR}, a far-going extension  of tilt stability and possibly improve results developed recently in \cite{MN2,MNR}. Furthermore, due to the strict connection of tilt stability and full stability to {\em strong stability} in the sense of Kojima \cite{ko}, which is equivalent of SSOSC under MFCQ as discussed in \cite[Chapter 5]{BS}, studying strong stability under weaker conditions than MFCQ, e.g., MSCQ (see also our Corollary~\ref{ssosc}) will be an interesting topic that we will pursue.

\end{document}